\documentclass{amsart}
\usepackage{etex}

\usepackage{amsmath,amsthm,amssymb,amscd,mathtools}
\usepackage{amsfonts}
\usepackage{rotating}
\usepackage{euscript}
\usepackage{pst-node}
\usepackage{epsfig,verbatim}
\usepackage{pictexwd,dcpic}
\usepackage{enumerate}
\usepackage{caption}
\usepackage{verbatim}

\def\be{\begin{equation}}
\def\ee{\end{equation}}

\def\C{{\mathbb C}} 
\def\f{\EuScript}
\def\N{{\mathbb N}} 
\def\P{{\mathbb P}}

\def\e{\eqref}
\def\phi{{\varphi}}
\def\v{{\varepsilon}} 
\def\tt{\widetilde}
\def\deg{{\rm deg\,}}

\def\Ker{{\rm Ker\,}}

\def\cos{{\rm cos\,}} 
 
\def\GCD{{\rm GCD }}

\def\Aut{{\rm Aut }}

\def\bp{\begin{proposition}}
\def\ep{\end{proposition}}

\def\bt{\begin{theorem}}
\def\et{\end{theorem}}
\def\br{\begin{remark}}
\def\er{\end{remark}}
\def\be{\begin{equation}}
\def\bee{\begin{equation*}}
\def\l{\label}

\def\e{\eqref}

\def\ee{\end{equation}}
\def\eee{\end{equation*}}
\def\bl{\begin{lemma}}
\def\el{\end{lemma}}
\def\bc{\begin{corollary}}
\def\ec{\end{corollary}}
\def\pr{\noindent{\it Proof. }}

\def\bd{\begin{definition}}
\def\ed{\end{definition}}
\def\t{\widetilde}

\def\hat{\widehat}

\newtheorem{theorem}{Theorem}[section]
\newtheorem{lemma}[theorem]{Lemma}
\newtheorem{corollary}[theorem]{Corollary}
\newtheorem{proposition}[theorem]{Proposition}
\theoremstyle{definition}
\newtheorem{remark}[theorem]{Remark}

\mathtoolsset{showonlyrefs}

\begin{document}
\title[Finiteness theorems for 
semiconjugate rational functions]{Finiteness theorems for commuting and
semiconjugate rational functions}
\author{Fedor Pakovich}
\thanks{This research was partially supported by the ISF, Grants No. 1432/18}

\begin{abstract}
Let $B$ be a fixed  rational function of one complex variable of degree at least two. In this paper, we study solutions of the functional equation $A\circ X=X\circ B$ in rational functions  $A$ and $X$. Our main result states that, unless $B$ is  a Latt\`es map or is 
conjugate to $z^{\pm d}$ or $\pm T_d$, the set of solutions  is finite, up to some natural transformations. 
In more detail, we show that there exist   finitely many rational functions $A_1, A_2,\dots, A_r$ and  $X_1, X_2,\dots, X_r$   such that the equality $A\circ X=X\circ B$ holds if and only if there exists a M\"obius transformation $\mu$ such that 
$A=\mu \circ A_j\circ \mu^{-1}$ and  $X=\mu \circ X_j\circ B^{\circ k}$  for some  $j,$ $1\leq j \leq r,$ 
and $k\geq 1$.  We also show that the number $r$ and the degrees $\deg X_j,$ $1\leq j \leq r,$  can be bounded from above in terms of the degree of $B$ only.
As an application, we prove an effective version of the  classical theorem of Ritt  about commuting rational functions.

\end{abstract}

\maketitle

\begin{section}{Introduction}
Commuting rational functions of one complex variable,  that is, rational solutions of the functional equation 
\be \l{comm} B\circ X=X\circ B\ee 
were investigated already at the dawn  of complex dynamics 
in the papers of Fatou, Julia, and Ritt 
\cite{f}, \cite{j}, \cite{r}. 
The 
most general result was obtained by Ritt. Roughly speaking, it  
states  that solutions of equation \eqref{comm} having no iterate in common
reduce either to powers, or to Chebyshev polynomials,  
or to Latt\`es maps.  
More precisely, in its modern formulation  due to Eremenko (\cite{e2}),
the Ritt theorem asserts that if $X$ and $B$ are commuting rational functions of degree at least two having no iterate in common, then there exists 
an {\it orbifold} $\f O$ of zero Euler characteristic defined on $\C\P^1$, $\C$, or $\C^*$ such that $A:\f O\rightarrow \f O$ and $B:\f O\rightarrow \f O$ are {\it covering maps} between orbifolds.
Notice that the Ritt theorem 
provides no information about commuting rational functions such that  
\be \l{ins}  X^{\circ l}= B^{\circ k}\ee for some $l,k\geq 1,$ and 
a characterization  of pairs $X$ and $B$ (commuting or not)  
satisfying \eqref{ins} is known only in the  polynomial case (\cite{r3}, \cite{r}).
Simple examples of commuting rational functions $X$ and $B$ satisfying \eqref{ins} can be obtained setting 
$$X=\mu_1 \circ R^{\circ l_1}, \ \ \ \ B=\mu_2 \circ R^{\circ l_2},$$ where $R$ is an arbitrary rational function and 
$\mu_1$, $\mu_2$ are M\"obius transformations commuting with $R$ and between themselves. However, it was shown already by Ritt (\cite{r}) that 
other examples also exist.

Functional equation \eqref{comm} 
is a particular case of the functional equation   
\be \l{i1}
A\circ X=X\circ B, 
\ee where $A$, $B$, and $X$ are rational functions, playing along with equation \eqref{comm} an important role in complex and arithmetic dynamics (see e.g. \cite{bu}, \cite{e2},  \cite{e3}, \cite{i}, \cite{ms}, \cite{pj}, \cite{aol}). 
We always will assume that  $A$ and $B$ in \eqref{i1} have degree at least two,  while $X$ has degree at least one.
In case if \eqref {i1} is satisfied for some  $X$ with $\deg X\geq 2$, 
the function $B$ is called {\it semiconjugate} to the function $A$, and the function $X$ is called a {\it semiconjugacy} from 
$B$ to $A.$ The case $\deg X=1$ corresponds to   the usual conjugacy. 
 In terms of dynamical systems, the conjugacy condition means that the  dynamical 
systems $B^{\circ k},$ $k\geq 1,$ and $A^{\circ k},$ $k\geq 1,$ are {\it equivalent}, 
while the semiconjugacy condition means that the system $A^{\circ k},$ $k\geq 1,$ is a {\it factor} of the
system $B^{\circ k},$ $k\geq 1.$

Semiconjugate  rational functions were investigated at length in the recent  papers \cite{semi},  \cite{arn}, \cite{lattes}. In particular, it was shown in \cite{semi} that 
if a solution $A,B,X$ of equation \eqref {i1} is  {\it primitive}, 
that is, satisfying the condition   
$\C(X,B)=\C(z)$, 
then there exist orbifolds $\f O_1$ and $\f O_2$ of non-negative Euler characteristic defined on $\C\P^1$ 
such that $A:\f O_1\rightarrow \f O_1$, $B:\f O_2\rightarrow \f O_2$, and $X:\f O_1\rightarrow \f O_2$ are {\it minimal holomorphic maps} between orbifolds. 
This condition generalizes the condition provided by the Ritt theorem, and implies strong restrictions on the possible  form of $A$, $B$ and $X$. In particular, it implies that the Galois closure of the field extension $\C(z)/\C(X)$ has 
genus zero or one.

Any solution of \eqref {i1}  reduces to a primitive one by a simple iterative process. 
Indeed, if  $\C(X,B)\neq \C(z)$, then 
by the L\"uroth theorem $\C(X,B)=\C(U_1)$
for some rational function $U_1$ of degree greater than one, and hence
\be \l{of} X=X_1\circ U_1,\ \ \  B=V_1\circ U_1  \ee 
for some rational functions $X_1$ and $V_1$. 
Substituting these expressions in \eqref{i1} we see  that the triple $A, X_1,U_1\circ V_1$ is another solution of \e{i1}.  This new
solution is not necessary primitive. Nevertheless, $\deg X_1 < \deg  X.$ Therefore,
after a finite number of similar transformations we will arrive at a primitive
solution.

In this paper, we study the totality of solutions of \eqref{i1} in rational functions $A$ and $X$  for a {\it fixed} rational function $B$. We 
show that, unless $B$ has a very special form,  the number of solutions of \eqref{i1}, considered up to some natural transformations, 
is {\it finite}. Moreover, this number can be bounded from above in terms of the degree of $B$ only. 
Notice that the results of  \cite{semi},  \cite{arn}, \cite{lattes}
do not immediately imply any  results of this kind,  since a priori the number of  steps 
in the reduction to a primitive solution as well as the number of primitive solutions for a  fixed function $B$ can be arbitrarily large.

In more details, let $A$, $X$ be a solution of \eqref{i1}.  
It is clear that $A$  is defined by $X$ in a unique way. However, for a given rational function $A$ there might be several $X$ satisfying \eqref{i1}.
In particular, for any $k\geq 1$ we may replace $X$ by  $\t X=X\circ B^{\circ k}.$
More generally,  if  $A,X$ is a solution  of \eqref {i1}, then for any M\"obius transformation $\mu$ and $k\geq 1$  we obtain another solution setting
\be \l{ebs}\t A= \mu \circ A\circ \mu^{-1}, \ \ \ \t X=\mu \circ X\circ B^{\circ k}.\ee
We say that a rational function $B$ is {\it special} if it is either  a Latt\`es map, or it is  conjugate to $z^{\pm n}$ or  $\pm T_n$.
In these terms,
the main  result of the paper is following.

\bt \l{main1} 
Let $B$ be a  non-special rational function of degree at least two.
Then there exist rational functions $A_1, A_2,\dots, A_{r}$ and  $X_1, X_2,\dots, X_{r}$
such that rational functions $A$ and $X$ satisfy \eqref {i1} if and only if 
$$ A=\mu \circ A_j\circ \mu^{-1}, \ \ \  X=\mu \circ X_j\circ B^{\circ k}$$
for some  $j,$ $1\leq j \leq r,$ M\"obius transformation $\mu$, and  $k\geq 0.$
 Furthermore, there exist (computable) functions $\psi,\phi:\N\rightarrow \N$ such that 
for any non-special $B$ of degree $d$ the number 
$r$ and the degrees $\deg X_j,$ $1\leq j \leq r,$ are bounded  from above by 
$\psi(d)$ and  $\phi(d),$ correspondingly. 
\et

Notice that Theorem \ref{main1} bounds both the number of conjugacy classes of rational functions $A$ such that \eqref{i1} holds  for some $X$, and, up to the second transformation in \eqref{ebs}, the number of rational functions $X$ such that \eqref {i1} holds for a fixed rational function $A$. 
In particular, for equation \eqref{comm} considered as a particular case of equation \eqref {i1} Theo\-rem \ref{main1} implies the following result.

\bt \l{main2} 
Let $B$ be a non-special rational function of degree at least two. Then there exist rational functions $X_1, X_2,\dots, X_{s}$  commuting with $B$ such 
that a rational function $X$ commutes with $B$ if and only if \be \l{comm+} X= X_j\circ B^{\circ k}\ee for some  $j,$ $1\leq j \leq s,$ and  $k\geq 0.$ 
Furthermore, there exist (computable) functions $\delta,\nu:\N\rightarrow \N$ such that 
for any non-special $B$ of degree $d$ the number 
$s$ and the degrees $\deg X_j,$ $1\leq j \leq r,$ are bounded from above by 
$\delta(d)$ and  $\nu(d),$ correspondingly. 
\et

Notice that Theorem \ref{main2} immediately implies  the Ritt theorem
about commuting rational functions in its part concerning non-special functions. Indeed, if $X$ commutes with  $B$, then 
any iterate $X^{\circ l}$ does. Thus, it follows from Theorem \ref{main2} by the Dirichlet box principle 
that there exist  distinct $l_1,$ $l_2\leq \delta(d)+1$ such that 
$$X^{\circ l_1}=X_j\circ B^{\circ k_1}, \ \ \ \ \ \ \ X^{\circ l_2}=X_j\circ B^{\circ k_2}$$ for the  same $j$ and 
some  $k_1,$ $k_2\geq 0$.  Assuming that $l_2>l_1,$ this yields 
\be \l{yuy} X^{\circ l_2}=X^{\circ l_1}\circ  B^{\circ (k_2-k_1)},\ee
implying that equality \eqref{ins} holds for $l=l_2-l_1$ and $k=k_2-k_1,$ since $X$ and $B$ commute. Moreover, 
  $l$ satisfies the inequality $l\leq \delta(d)$. Thus, Theorem \ref{main2} improves the Ritt theorem, which provides 
neither existence of finitely many functions  such that any  function $X$ commuting with $B$	has form \eqref{comm+}, nor boundedness of $l$ in terms of $\deg B.$ 


The paper is organized as follows. In the second section,  
we fix the notation and  recall main definitions and results related to Riemann surface orbifolds. We also collect some technical results, mostly from the papers \cite{semi}, \cite{lattes}, used in the following. In the third section, we 
study the systems of functional equations
\be \l{asd0} U_{i}\circ V_{i}=V_{i+1}\circ U_{i+1},\ \ \ 1\leq i \leq s-1,\ee 
where  $U_i,$ $V_i,$ $1\leq i \leq s,$ are rational functions of degree at least two, which we will call {\it chains}.
Such chains correspond to chains 
 \be \l{chh} B\rightarrow B_1 \rightarrow B_2  \rightarrow \dots \rightarrow B_s\ee 
of rational functions
$$B=V_1\circ U_1, \ \ \  B_i= U_i\circ V_i, \ \ \ \ \ 1\leq i\leq s,$$ 
relating an arbitrary  solution of \eqref{i1} with a primitive one,
and the main result of the third section asserts that under certain restrictions the length $s$ of such a chain can be effectively bounded in terms of the degree of $B$.

In the fourth section, we define  
{\it an extended symmetry group} of a rational function $F$
as the group of M\"obius transformations $\sigma$ such that 
$$ F\circ \sigma=\nu \circ F$$ for some  M\"obius transformations $\nu$. We show that,  unless $F=\mu_1\circ z^n\circ \mu_2$ for some M\"obius transformations $\mu_1,$ $\mu_2$, this group is finite, and using this fact prove Theorem \ref{main1} for primitive solutions of \eqref{i1}. 
In the fifth section, we prove a ``compositional'' counterpart of Theorem \ref{main1}, which asserts that if $B$ is  non-special, then any $X$  such that \eqref{i1} holds for some $A$ 
can be decomposed  into a composition of rational functions  $$X=X^{\prime}\circ U\circ B^{\circ k},\ \ \ k\geq 1,$$ such that the Galois closure of $\C(z)/\C(X')$ has 
genus zero or one, and  $U$ is a ``compositional right factor'' of some iterate $B^{\circ l}$  
with $l$ bounded in terms of $\deg B.$ We also prove  Theorem \ref{main1} and Theorem  \ref{main2}. Finally, we prove an effective version of the Ritt theorem.
\end{section}

\begin{section}{Functional decompositions and orbifolds}

\begin{subsection}{Orbifolds and maps between orbifolds}
In this section we fix the notation and  recall main definitions and results related to Riemann surface orbifolds  (see  \cite{mil}, Appendix E). We also collect some technical results,  from the papers  \cite{semi}, \cite{gen}, \cite{gen0}, \cite{lattes}, used in the following. 

A pair $\f O=(R,\nu)$ consisting of a Riemann surface $R$ and a ramification function $\nu:R\rightarrow \mathbb N$ which takes the value $\nu(z)=1$ except at isolated points is called an {\it orbifold}. For an orbifold $\f O$ the {\it  Euler characteristic} of $\f O$ is the number
\be \l{echa} \chi(\f O)=\chi(R)+\sum_{z\in \C\P^1}\left(\frac{1}{\nu(z)}-1\right),\ee
the set of {\it singular points} of $\f O$ is the set 
$$c(\f O)=\{z_1,z_2, \dots, z_s, \dots \}=\{z\in \C\P^1 \mid \nu(z)>1\},$$ and  the {\it signature} of $\f O$ is the set 
$$\nu(\f O)=\{\nu(z_1),\nu(z_2), \dots , \nu(z_s), \dots \}.$$

If $R_1$, $R_2$ are Riemann surfaces provided with ramification functions $\nu_1,$ $\nu_2$, and 
$f:\, R_1\rightarrow R_2$ is a holomorphic branched covering map, then $f$
is called  {\it a covering map} $f:\,  \f O_1\rightarrow \f O_2$
between orbifolds
$\f O_1=(R_1,\nu_1)$ and $\f O_2=(R_2,\nu_2)$
if for any $z\in R_1$ the equality 
\be \l{us} \nu_{2}(f(z))=\nu_ {1}(z)\deg_zf\ee holds, where $\deg_zf$ stands for the local degree of $f$ at the point $z$.
If for any $z\in R_1$ instead of equality \eqref{us} 
the weaker condition 
\be \l{uuss} \nu_{2}(f(z))\mid \nu_ {1}(z)\deg_zf\ee
holds,  then the map $f$
is called {\it a holomorphic map} $f:\,  \f O_1\rightarrow \f O_2$
between orbifolds
$\f O_1=(R_1,\nu_1)$ and $\f O_2=(R_2,\nu_2).$

A universal covering of an orbifold ${\f O}=(R,\nu)$
is a covering map between orbifolds  $\theta_{\f O}:\,
\tt {\f O}\rightarrow \f O$ such that $\tt R$ is simply connected and $\tt \nu(z)\equiv 1.$ 
If $\theta_{\f O}$ is such a map, then 
there exists a group $\Gamma_{\f O}$ of conformal automorphisms of $\tt R$ such that the equality 
$\theta_{\f O}(z_1)=\theta_{\f O}(z_2)$ holds for $z_1,z_2\in \tt R$ if and only if $z_1=\sigma(z_2)$ for some $\sigma\in \Gamma_{\f O}.$ A universal covering exists and 
is unique up to a conformal isomorphism of $\tt R,$
unless $\f O$ is the Riemann sphere with one ramified point, or  the Riemann sphere with two ramified points $z_1,$ $z_2$ such that $\nu(z_1)\neq \nu(z_2)$  (see \cite{fk}, 
Section IV.9.12).
Abusing  notation we will denote by $\tt {\f O}$ both the
orbifold and the  Riemann surface  $\tt R$.

Covering maps between orbifolds lift to isomorphisms between their universal coverings.
More generally, for holomorphic maps between orbifolds the following proposition holds  (see \cite{semi}, Propo\-si\-tion 3.1).

\bp \l{poiu} Let $f:\,  \f O_1\rightarrow \f O_2$ be a holomorphic map between orbifolds. Then for any choice of $\theta_{\f O_1}$ and $\theta_{\f O_2}$ there exist 
a holomorphic map $F:\, \tt {\f O_1} \rightarrow \tt {\f O_2}$ and 
a homomorphism $\phi:\, \Gamma_{\f O_1}\rightarrow \Gamma_{\f O_2}$ such that the diagram 
\be \l{dia2}
\begin{CD}
\tt {\f O_1} @>F>> \tt {\f O_2}\\
@VV\theta_{\f O_1}V @VV\theta_{\f O_2}V\\ 
\f O_1 @>f >> \f O_2\ 
\end{CD}
\ee
is commutative and 
for any $\sigma\in \Gamma_{\f O_1}$ the equality
\be \l{homm}  F\circ\sigma=\phi(\sigma)\circ F \ee holds.
The map $F$ is defined by $\theta_{\f O_1}$, $\theta_{\f O_2}$, and $f$  
uniquely up to a transformation 
$F\rightarrow g\circ F,$ where $g\in \Gamma_{\f O_2}$. 
In the other direction, for any holomorphic map  $F:\, \tt {\f O_1} \rightarrow \tt {\f O_2}$  which satisfies \eqref{homm} for some homomorphism $\phi:\, \Gamma_{\f O_1}\rightarrow \Gamma_{\f O_2}$
there exists a uniquely defined  holomorphic map between orbifolds $f:\,  \f O_1\rightarrow \f O_2$ such that diagram \eqref{dia2} is commutative.
The holomorphic map $F$ is an isomorphism if and only if $f$ is a covering map between orbifolds. \qed

\ep

If $f:\,  \f O_1\rightarrow \f O_2$ is a covering map between orbifolds  $\f O_1$ and $\f O_2$ with compact supports $R_1$ and $R_2$, 
then  the Riemann-Hurwitz 
formula implies that 
\be \l{rhor} \chi(\f O_1)=d \chi(\f O_2), \ee
where $d=\deg f$ (see \cite{mil}). 
For holomorphic maps the following statement is true (see \cite{semi}, Proposition 3.2). 

\bp \l{p1} Let $f:\, \f O_1\rightarrow \f O_2$ be a holomorphic map between orbifolds with compact supports.
Then 
\be \l{iioopp} \chi(\f O_1)\leq \chi(\f O_2)\,\deg f, \ee and the equality 
holds if and only if $f:\, \f O_1\rightarrow \f O_2$ is a covering map between orbifolds. \qed
\ep

Let $R_1$, $R_2$ be Riemann surfaces and 
$f:\, R_1\rightarrow R_2$ a holomorphic branched covering map. Assume that $R_2$ is provided with ramification function $\nu_2$. In order to define a ramification function $\nu_1$ on $R_1$ so that $f$ would be a holomorphic map between orbifolds $\f O_1=(R_1,\nu_1)$ and $\f O_2=(R_2,\nu_2)$ 
we must satisfy condition \eqref{uuss}, and it is easy to see that
for any  $z\in R_1$ a minimum possible value for $\nu_1(z)$ is defined by 
the equality 
\be \l{rys} \nu_{2}(f(z))=\nu_ {1}(z)\GCD(\deg_zf, \nu_{2}(f(z)).\ee 
In case if \eqref{rys} is satisfied for  any $z\in R_1$ we 
say that $f$ is {\it a  minimal holomorphic  map} 
between orbifolds 
$\f O_1=(R_1,\nu_1)$ and $\f O_2=(R_2,\nu_2)$.

It follows from the definition that for any orbifold $\f O=(R,\nu)$ and a holomorphic branched covering map $f:\, R^{\prime} \rightarrow R$ there exists a unique orbifold structure \linebreak $\f O^{\prime}=(R^{\prime},\nu^{\prime})$  such that 
$f:\f O^{\prime}\rightarrow \f O$ is a minimal holomorphic map between orbifolds. 
We will denote the corresponding orbifold by $f^*\f O.$ Notice that any covering map between orbifolds $f:\,  \f O_1\rightarrow \f O_2$ is a  minimal holomorphic map. In particular, this implies that for any covering map  $f:\,  \f O_1\rightarrow \f O_2$ the equality
$\f O_1=f^*\f O_2$ holds.
Minimal holomorphic maps between orbifolds possess the following fundamental property with respect to the operation of  composition (see \cite{semi}, Theorem 4.1).

\bt \l{serrr} Let $f:\, R^{\prime\prime} \rightarrow R^{\prime}$ and $g:\, R^{\prime} \rightarrow R$ be holomorphic branched covering maps, and  $\f O=(R,\nu)$ an orbifold. 
Then 

$$(g\circ f)^*\f O= f^*(g^*\f O).\eqno{\Box}$$
\et

Theorem \ref{serrr} implies in particular the following corollaries (see   \cite{semi}, Corollary 4.1 and Corollary 4.2).

\bc \l{serka0} Let $f:\, \f O_1\rightarrow \f O^{\prime}$ and $g:\, \f O^{\prime}\rightarrow \f O_2$ be minimal holomorphic maps (resp. covering maps) between orbifolds.
Then  $g\circ f:\, \f O_1\rightarrow \f O_2$ is  a minimal holomorphic map (resp. covering map). \qed
\ec

\bc \l{indu2}  Let $f:\, R_1 \rightarrow R^{\prime}$ and $g:\, R^{\prime} \rightarrow R_2$ be holomorphic branched covering maps, and  $\f O_1=(R_1,\nu_1)$ and  $\f O_2=(R_2,\nu_2)$
orbifolds. Assume that \linebreak $g\circ f:\, \f O_1\rightarrow \f O_2$ is  a minimal holomorphic map (resp. a co\-vering map). Then $g:\, g^*\f O_2\rightarrow \f O_2$  and  $f:\, \f O_1\rightarrow g^*\f O_2 $ are minimal holomorphic maps (resp. covering maps). \qed
\ec

In this paper, essentially all considered orbifolds will be defined on $\C\P^1.$ 
So, we will omit  the Riemann surface $R$ in the definition of $\f O=(R,\nu)$
meaning that $R=\C\P^1.$  
``Most'' orbifolds on  $\C\P^1$ have negative Euler characteristic. Orbifolds $\f O$  with $\chi(\f O)\geq 0$  
and corresponding $\Gamma_{\f O}$ and $\theta_{\f O}$ can be described explicitly as follows. 
The equality $\chi(\f O)=0$ holds if and only if the signature of $\f O$ 
belongs to the list
\be \l{list}\{2,2,2,2\} \ \ \ \{3,3,3\}, \ \ \  \{2,4,4\}, \ \ \  \{2,3,6\}, \ee while $\chi(\f O)>0$  if and only if either $\f O$ is the non-ramified sphere or
 the signature of $\f O$  belongs to the  list 
 \be \l{list2} \{n,n\}, \ \ n\geq 2,  \ \ \ \{2,2,n\}, \ \ n\geq 2,  \ \ \ \{2,3,3\}, \ \ \ \{2,3,4\}, \ \ \ \{2,3,5\}.\ee

Groups $\Gamma_{\f O}\subset \Aut(\C)$ corresponding to orbifolds $\f O$ with signatures \eqref{list}  
are generated by translations of $\C$ by elements of some lattice $L\subset \C$ of rank two and the rotation $z\rightarrow  \v z,$ where $\v$ is an $n$th root of unity with $n$ equal to 2,3,4, or 6, such that  $\v L=L$.  
 Accordingly, the functions $\theta_{\f O}$ 
may be written in terms of the  corresponding
Weierstrass functions as $\wp(z),$ $\wp^{\prime }(z),$ $\wp^2(z),$  and $\wp^{\prime 2}(z)$  (see  \cite{mil2}, or \cite{fk}, 
Section IV.9.5).
Groups $\Gamma_{\f O}\subset \Aut(\C\P^1)$ corresponding to   orbifolds $\f O$ with signatures \eqref{list2} are the well-known five finite subgroups 
 $C_n,$  $D_{2n},$  $A_4,$ $S_4,$ $A_5$ of $\Aut(\C\P^1)$, and the functions $\theta_{\f O}$ are Galois coverings of $\C\P^1$ by $\C\P^1$ of degrees 
$n$, $2n,$ $12,$ $24,$  $60,$ calculated for the first time by Klein in \cite{klein}.
In particular, for $C_n$ and $D_{2n}$ the corresponding functions $\theta_{\f O}$ are $z^n$ and 
 $$Z_n= \frac {1}{2}\left(z^n+\frac {1}{z^n}\right).$$

\end{subsection}

\begin{subsection}{Orbifolds $\f O_1^A$, 
$\f O_2^A$, and  $\f O_0^A$.}

With each rational  function $A$ 
one can associate in a natural way two orbifolds $\f O_1^A$
and 
$\f O_2^A$
setting $\nu_2^A(z)$  
equal to the least common multiple of local degrees of $f$ at the points 
of the preimage $A^{-1}\{z\}$, and $$\nu_1^A(z)=\nu_2^A(A(z))/\deg_zA.$$
By construction,  $A:\, \f O_1^A\rightarrow \f O_2^A$ 
is a covering map between orbifolds.  Furthermore, since 
the composition $A\circ \theta_{\f O_1^A}: \t{\f O_1^A}\rightarrow \f O_2^A$ is a covering map  between orbifolds by Corollary \ref{serka0}, it follows from the uniqueness of the universal co\-vering that 
\be \l{ravv}  \theta_{\f O_2^A}=A\circ \theta_{\f O_1^A}.\ee

We recall that 
a {\it Latt\`es map} can be defined as a rational function $A$ of degree at least two such that  $A:\f O\rightarrow \f O$ is a  covering  map
for some orbifold $\f O$  (see \cite{mil2}, \cite{lattes}). Such an orbifold is defined in a unique way and necessarily 
satisfies the condition $\chi(\f O)=0$ in view of equality \eqref{rhor}.   
Following \cite{lattes}, we say that  a rational function $A$ of degree at least two is a {\it generalized Latt\`es map} if there exists an orbifold $\f O$ distinct from the non-ramified sphere such that  $A:\,  \f O\rightarrow \f O$ is a  minimal holomorphic map between orbifolds.
Notice that, similarly to usual Latt\`es maps, generalized Latt\`es maps can be 
described in terms of group actions and semiconjugacies (see \cite{lattes}).

In general, there might be more than one orbifold $\f O$ such that $A:\, \f O \rightarrow \f O$ is a minimal holomorphic map  between orbifolds, and even infinitely many such orbifolds.  
Namely, the power $z^{\pm d}:\f O\rightarrow \f O$ is a minimal holomorphic map for any $\f O$ defined by the conditions \be \l{raz} 
\nu(0)=\nu(\infty)=n, \ \ \ \ n\geq 2, \ \ \ \  \GCD(d,n)=1,\ee and the Chebyshev polynomial $\pm T_{d}:\f O\rightarrow \f O$ of degree $d$ is a minimal holomorphic map  for any $\f O$ defined by  the conditions \be \l{dva} 
\nu(-1)=\nu(1)=2, \ \ \ \nu(\infty)=n,\ \ \ \ n\geq 1, \ \ \ \  \GCD(d,n)=1.\ee 
Nevertheless, if $A$ is  not conjugate to $z^{\pm n}$ or $\pm T_n$, there exists  
a ``maximal'' orbifold $\f O$, denoted by $\f O_0^A$, such that $A:\,  \f O\rightarrow \f O$ is a  minimal holomorphic map. 

In more details, 
for   orbifolds $\f O_1$  and $\f O_2$ we 
write $ \f O_1\preceq \f O_2$ if for any $z\in\C\P^1$
the condition   $\nu_1(z)\mid \nu_2(z)$ holds.
In this notation, the following statement holds (see \cite{lattes}, Theorem 1.2).

\bt \l{uni} 
Let $A$ be a rational function  of degree at least two not conjugate to $z^{\pm d}$ or $\pm T_d.$ Then there exists an orbifold $\f O_0^A$ such that $A:\, \f O_0^A\rightarrow \f O_0^A$
is a minimal holomorphic map between orbifolds, and for any orbifold $\f O$ such that 
$A:\, \f O\rightarrow \f O$ is a minimal holomorphic map between orbifolds the relation $\f O\preceq \f O_0^A$ holds. Furthermore, $\f O_0^{A^{\circ l}}=\f O_0^A$ for any 
$l\geq 1$.
\et

Clearly, generalized Latt\`es maps are exactly rational functions for which the orbifold $\f O_0^A$ is distinct from the non-ramified sphere, completed
by the functions  $z^{\pm d}$ or $\pm T_d$ for which the orbifold $\f O_0^A$ is not defined. Furthermore, a rational function $A$ is a  Latt\`es map if and only if  $\chi(\f O_0^A)= 0$ (see \cite{lattes}, Lemma 6.4).

\end{subsection}

\begin{subsection}{Functions $A$ with $\chi(\f O_2^A)\geq 0$}\label{kot}  

Rational functions $A$ for which the orbifold $\f O_2^A$ has non-negative Euler characteristic play a special role in the description of solutions of \eqref{i1}. Below we list some properties of such functions used below. 

Let $F$ and $G$ be rational functions. We will
call $G$ a {\it compositional left factor} of $F$ if $F = G\circ H$ for some rational function $H$. Compositional right factors are defined in a similar way. 
We will say that rational functions $A_1$ and $A_2$ are $\mu$-equivalent, and write 
$${A_1}\underset{\mu}{\sim} A_2,$$ if $A_1$ and $A_2$  satisfy the equality
$$A_1= \mu_1\circ A_2\circ \mu_2,$$
for some M\"obius transformations
$\mu_1$ and $\mu_2$.

Recall that for a rational function $X$  
its  normalization $\t X$ is defined as a holomorphic function of the lowest possible degree
between compact Riemann surfaces  $\t X:\,\t S_X\rightarrow \C\P^1$  such that $\t X$ is a Galois covering and
 $$\t X=X\circ H$$ for some  holomorphic map $H:\,\t S_X\rightarrow \C\P^1$. From the algebraic point of view the passage from $X$ to $\t X$ 
corresponds to the passage from the field extension $\C(z)/\C(X)$ to its Galois closure. 

In the above  terms, rational functions $A$ for which $\chi(\f O_2^A)\geq 0$ can be characterized as follows  (see \cite{gen0}, Lemma 2.1).

\bl \l{ml} 
Let $A$ be a  rational function. Then $g(\t S_A)=0$ if and only if \linebreak $\chi(\f O_2^A)> 0$, and  $g(\t S_A)=1$ if and only if  $\chi(\f O_2^A)= 0$. \qed
\el 
Since, by Lemma \ref{ml}, rational functions $A$ with $\chi(\f O_2^A)> 0$ are 
compositional left factors of rational Galois coverings of  $\C\P^1$ by $\C\P^1$, these functions  can be listed explicitly (see Theorem 1.1 in \cite{gen0}). 
Below we will  need only the following corollary of this classification.

\bl \l{xxx} Let $A$ be a rational function of degree at least two such that \linebreak $\chi(\f O_2^A)>0.$ 
Then either 
$\deg A$ belongs to the set 
$$E_1=\{4,5,6,8,10,12,15,20,24,30,60\},$$ or $A$ is $\mu$-equivalent to one of the functions $z^{n}$,  $T_{n}$,  $Z_{n}$, where $n\geq 2.$ \qed
\el

Notice that 
$$z^2 \underset{\mu}{\sim} T_2\underset{\mu}{\sim}Z_1.$$
However,  for
$n>2$ the functions $z^n$, $T_n,$ and $Z_{n/2}$ are   pairwise not $\mu$-equivalent.

\vskip 0.2cm
Consider now rational functions $A$ with $\chi(\f O_2^A)=0.$ Since $A:\f O_1^A\rightarrow \f O_2^A$ is a covering map, 
and $\chi(\f O_1^A)=0$ by \eqref{rhor}, any such a map 
 is a covering map between orbifolds of zero Euler characteristic. In the other direction, it can be  shown 
that  if $A\,: \f O_1\rightarrow \f O_2$   is a covering  map between 
orbifolds of zero  Euler characteristic, then    with a few exceptions
the equalities $$\f O_1^A=\f O_1, \ \ \ \f O_2^A=\f O_2$$ hold (see \cite{gen0}, Theorem 5.2).  Again, we will need only the following corollary of this result.

\bl \l{xxxx} Let $A$ be a rational function of degree at least two, and  let 
 $\f O_1$, $\f O_2$ be orbifolds such that $A:\, \f O_1\rightarrow \f O_2$ is a covering map between orbifolds and
$\chi(\f O_1)=0,$ $\chi(\f O_2)=0.$ Then either $\deg A$ belongs to the set 
$$E_2=\{2,3,4,6,8,12\},$$ or the equalities 
$\f O_2=\f O_2^A$, $\f O_1=\f O_1^A$ hold. \qed
\el

\vskip 0.2cm

Finally, we will use the following well-known statement concerning decompositions of the functions $z^{n}$,  $T_{n}$,  $Z_{n}$
into a composition of two rational functions (see e.g. \cite{gen0}, Sections 4.1-4.2).

\bl \l{any} \ \ 
\vskip 0.1cm
\begin{enumerate}[a)]

\item Any decomposition of $z^n$, $n\geq 2,$   has the form \be \l{ega} z^n=(z^{n/d}\circ \mu)\circ(\mu^{-1}\circ  z^d),\ee where $d\vert n$ and $\mu$ is a M\"obius transformation.

\item Any decomposition of  $Z_n$, $n\geq 2,$  either has the form 
\be \l{ega1}Z_n=\left(Z_{n/d}\circ \mu\right)\circ\left(\mu^{-1}\circ z^d\right),\ee
where $d\vert n$ and  $\mu$ is a M\"obius transformation,
or has the form \be \l{ega2}
 Z_n=\left(\v^{n} T_{n/d}\circ \mu\right)\circ\left(\mu^{-1}\circ Z_d\circ (\v z)\right),
\ee where $d\vert n$, $\v^{2n}=1,$ and $\mu$ is a M\"obius transformation. 

\item Any decomposition of $T_n$, $n\geq 2,$  has the form \be \l{xerxyi} T_n=(T_{n/d}\circ \mu)\circ(\mu^{-1}\circ  T_d),\ee where $d\vert n$ and $\mu$ is a M\"obius transformation. \qed
\end{enumerate}
\el

\end{subsection}

\begin{subsection}{Equivalence $\sim$ and special functions}
Let $B$ be a rational function. 
For any decomposition $B=V\circ U,$ where $U$ and $V$ are rational functions, the 
rational function $\t B=U\circ V$ is called an {\it elementary transformation} of $B$. We say that rational functions $B$ and $A$ are {\it equivalent} and write $A\sim B$ if there exists 
a chain of elementary transformations between $B$ and $A$
(this equivalence relation  should not be confused with the equivalence relation from Section 2.3 where the
subscript $\mu$ is used). For a rational function $B$ we will denote its equivalence class by $[B].$ Since for any M\"obius transformation $W$ the equality
$$B=(B\circ W)\circ W^{-1}$$ holds, each equivalence class $[B]$ is a union of conjugacy classes. 

Equivalent functions provide examples of semiconjugate functions.  
Indeed, since for $B$ and $\t B$ as above the equalities 
$$\t B\circ U=U\circ B, \ \ \ \ \  B\circ V=V\circ \t B$$ hold,  
 $B$ is semiconjugate to $\t B$, and   $\t B$ is semiconjugate to $B,$
implying inductively that if
$A\sim B$, then $A$ is semiconjugate to $B$, and   $B$ is semiconjugate to $A$.


A rational function $B$ of degree at least two is called {\it special} if it is either a Latt\`es map, or it is  conjugate to $z^{\pm n}$ or  $\pm T_n$. Special functions can be characterized as  {\it  finite quotients
of affine maps} in the following sense: a rational function $B$ is special if and only if 
 there is a lattice of rank one or two $\Lambda\subset \C$, an affine map $L=at+b$ from $\C/\Lambda$ to $\C/\Lambda$,  
and a finite-to-one holomorphic map 
$\theta: \C/\Lambda\rightarrow \C\P^1\setminus {\f E}_B$,  where $\f E_B$ stands for the set of exceptional values of $B$,   which satisfy the
semiconjugacy relation 
\be 
\begin{CD}
 \C/\Lambda @>L>> \C/\Lambda \\ @VV\theta V @VV\theta V\\ 
\C\P^1\setminus {\f E}_B @>B >> \C\P^1\setminus {\f E}_B 
\end{CD}
\ee
(see \cite{mil2}). 

Equivalently, special functions can be described as rational functions $B$ that  
are covering maps $B:\f O\rightarrow \f O$ between orbifolds for some $\f O=(R,\nu)$ with 
$R=\C\P^1\setminus {\f E}_B$. It follows from \eqref{rhor} that for such an orbifold $\f O$ the equality $\chi(\f O)=0$ holds, and 
\eqref{echa} implies easily that,  if ${\f E}_B\neq \emptyset$, then either ${\f E}_B$ contains two points  and $\nu\equiv 1$, or ${\f E}_B$ contains one point and $\nu(\f O)=\{2,2\}.$ 
Correspondingly, the map $z^{\pm d}: \f O\rightarrow \f O$ is a covering map for the non-ramified orbifold with $R=\C\P^1\setminus\{0,\infty\}$, while 
$\pm T_d: \f O\rightarrow \f O$ is a covering map for the orbifold defined on $ R=\C\P^1\setminus \{\infty\}$ by the condition $\nu(1)=2,$ $\nu(-1)=2$. The corresponding functions $\theta$ are $e^z$ and $\cos z$. 


\vskip 0.2cm

Below we collect several facts about special functions that we will need in the following. 

\bl \l{korova} Let $F$ be a special rational function and $\hat F\sim  F.$
Then $\hat F$ is special.  
\el 
\pr Assume that $F$ is conjugate to $z^{\pm n}$.  Then, by Lemma  \ref{any},   any elementary transformation of $F$ is conjugate to  $z^{\pm n}$, implying inductively  that any $\hat F\sim F$ 
 is conjugate to $z^{\pm n}$.  If $F$ is conjugate to $\pm T_n$, the proof is similar. Finally, if $F$ is a Latt\`es, then $\hat F$ is a Latt\`es map (see \cite{lattes}, Corollary 4.4). \qed

\vskip 0.2cm

The next statement follows from Corollary 4.7 and Lemma 6.3 in the paper \cite{lattes}.

\bl \l{baran} Let $A$ be a rational function of degree $d\geq 2$ such that some iterate $A^{\circ l}$, $l\geq 2,$ is 
special. Then $A$ is special.\qed
\el

Finally, we will need the following result (see \cite{mil2}, Corollary 4.3).

\bl \l{zaq} Let $A$ and $B$ be rational functions such that $A$ is semiconjugate to $B.$ Then $A$ is special if and only if $B$ is special.  \qed
\el

\end{subsection}

\begin{subsection}{Good solutions of $A\circ C=D\circ B$}
We say that 
a solution $A,C,D,B$ of the functional equation 
\be \l{m} A\circ C=D\circ B\ee
 in rational functions is {\it good} if 
the algebraic curve 
\be \l{ad}\f E_{A,D}: A(x)-D(y)=0\ee is irreducible
and $\C(C,B)=\C(z)$.
This definition is a particular case of the definition of good solutions of
\eqref{m} in holomorphic functions defined on compact Riemann surfaces (see \cite{semi}, Section 2). In particular, Theorem 2.1 in \cite{semi} implies the following.

\bl \l{good2} If a 
solution  $A,C,D,B$  of \e{m} is good, then $\deg C=\deg D.$ \qed
\el

Furthermore, the following statement holds (see \cite{semi}, Lemma 2.1).

\bl \l{good} A
solution  $A,C,D,B$  of \e{m} is good
 whenever any two of the following three conditions are satisfied:
\begin{itemize}
\item the curve $\f E_{A,D}$ is irreducible,
\item the equality  $\C(C,B)=\C(z)$ holds,
\item the equality  $\deg C=\deg D$ holds.  \qed
\end{itemize}
\el

The property of a solution $A,C,D,B$ of \eqref{m} to be good imposes strong restrictions on the ramification collections 
of the functions $A,C,D,B$, which are described by the following theorem 
(see \cite{semi}, Theorem 4.2).

\bt \l{t1} Let $A,C,D,B$ be a  good solution of \eqref{m}. Then
the commutative diagram 
\be \l{diag} 
\begin{CD}
\f O_1^C @>B>> \f O_1^D\\
@VV C V @VV D V\\ 
\f O_2^C @>A >> \f O_2^D\ 
\end{CD}
\ee
consists of minimal holomorphic  maps between orbifolds. \qed
\et

Since an irreducible algebraic curve $\f E_{A,D}$ has genus zero if and only if it can be parametrized 
by some rational functions $C$ and $B$ with $\C(C,B)=\C(z)$, describing good solutions of \eqref{m} mostly reduces to describing irreducible algebraic curves \eqref{ad} of genus zero. The following general result is proved in \cite{gen}.

\bt \l{m2} 
Let $A$ be a rational function of degree $n$ such that $\chi(\f O_2^A)<0.$
Then for any rational function 
$D$ of degree $m$ such that the curve $\f E_{A,D}$  is irreducible 
the inequality 
\be \l{ma} 
g(\f E_{A,D}) >\frac{m-84n+168}{168}
\ee holds. \qed
\et 

The practical meaning of Theorem \ref{m2} is that whenever $A,C,D,B$ is a good solution of \eqref{m} with 
$$\deg D\geq 84(\deg A-2),$$ the function $A$ necessarily satisfies the restrictive condition  
$\chi(\f O_2^A)\geq 0$ discussed in Section \ref{kot}.

The next result we will need states  that ``gluing together'' two commutative diagrams corresponding to good solutions of \eqref{m} we obtain again a good solution of \eqref{m} (see the diagram below) 
\be 
\begin{CD}
\C\P^1 @>B>> \C\P^1  @>W>> \C\P^1 \\
@VV C V @VV D V @VV {V} V \\ 
\C\P^1 @>A >> \C\P^1  @>U >> \C\P^1\,. 
\end{CD}
\ee

\bt \l{sum+} Assume that 
$A,C,D,B$ and $U,D,V,W$ are  good solutions of \eqref{m}. 
Then $U\circ A$, $C$, $V,$ $W\circ B$ is also  a good solution of \eqref{m}. 
\et 
\pr The  theorem is a particular case of Theorem 2.10 in \cite{aol}. For the reader convenience 
we provide a short independent proof.

Since Lemma \ref{good2} implies the equalities 
\be \l{zza} \deg C=\deg V, \ \ \ \deg (W\circ B)=\deg (A\circ U),\ee it 
follows from Lemma \ref{good} that it is enough to prove that 
the curve  \be \l{ka} (U\circ A)(x)-V(y)=0\ee is irreducible. Assume the inverse. Then 
$$t\rightarrow (C(t), (W\circ B)(t))$$ is a parametrization of some proper irreducible 
component $F(x,y)=0$ of \eqref{ka}, implying that 
\be \l{zazaz} C=X\circ R, \ \ \ \ W\circ B=Y\circ R,\ee for some rational functions $X,Y$ and $R$ such that 
\be \l{azz} \deg X=\deg_yF, \ \ \ \ \deg Y=\deg_xF.\ee 
Moreover, $\deg R>1$ since otherwise  equalities \eqref{zza}, \eqref{zazaz}, 
 and \eqref{azz}  imply that the curve $F(x,y)=0$ coincides with \eqref{ka}.

Since $U,D,V,W$ is a  good solution of \eqref{m}, it follows from the
equality $$(U\circ A)\circ X= U\circ (A\circ X)= V\circ Y$$ 
 that 
there exists  a rational function $T$ such that 
$$A\circ X=D\circ T, \ \ \ \ Y=W\circ T.$$
Similarly, the first of these equalities implies that
there exists a rational function  $S$ 
such that 
$$X=C\circ S, \ \ \ \ T =B\circ S.$$
Thus, 
$$X=C\circ S, \ \  \ \ Y=W\circ B\circ S,$$
implying by \eqref{zazaz} that $\deg R= 1.$ The contradiction obtained shows that  \eqref{ka} is irreducible. \qed

\end{subsection}

\begin{subsection}{Primitive solutions of $A\circ X=X\circ B$}
Recall that  a solution $A,X,B$ of equation  \eqref{i1} is called {\it primitive} if $\C(X,B)=\C(z)$. By Lemma \ref{good}, a solution $A,X,B$ of  \eqref{i1}  is primitive if and only if the corresponding solution $$A=A, \ \ \ C=X, \ \ \  D=X, \ \ \ B=B$$ of \eqref{m} is good. Primitive solutions are described as follows (see \cite{semi}, Theorem 6.1).

\bt \l{oip} Let $A,B,X$ be rational functions of degree at least two  such that $A\circ X=X\circ B$ and  $\C(B,X)=\C(z)$.   Then $\chi(\f O_1^X)\geq 0$, $\chi(\f O_2^X)\geq 0$, and 
the commutative diagram 
\be 
\begin{CD} \l{gooopa}
\f O_1^X @>B>> \f O_1^X\\
@VV X V @VV X V\\ 
\f O_2^X @>A >> \f O_2^X\ 
\end{CD}
\ee
consists of minimal holomorphic  maps between orbifolds. \qed
\et

The following statement (see \cite{semi}, Theorem 5.1) is a 
more precise version of Proposition \ref{poiu} for minimal holomorphic maps $A:\f O\rightarrow \f O$ with $\chi(\f O)>0$. 


\bt \l{las}  Let $A$ and $F$ be rational functions of degree at least two and  $\f O$ an orbifold with $\chi(\f O)>0$ such that 
$A:\, \f O \rightarrow \f O$ 
is a holomorphic map between orbifolds  and  the diagram 
\be \l{dia3}
\begin{CD}
\tt {\f O} @>F>> \tt {\f O}\\
@VV\theta_{\f O}V @VV\theta_{\f O}V\\ 
\f O @>A >> \f O\ 
\end{CD}
\ee
commutes.
Then the following conditions are equivalent.

\begin{itemize}
\item[(1)] The holomorphic map $A$ is a minimal holomorphic  map. 
\item[(2)]  The homomorphism $\phi:\, \Gamma_{\f O}\rightarrow \Gamma_{\f O}$ defined by the equality  
$$F\circ\sigma=\phi(\sigma)\circ F, \ \ \ \sigma\in \Gamma_{\f O},$$ is an automorphism of $\Gamma_{\f O}$.
\item[(3)] The triple $F,$ $A,$ $\theta_{\f O}$ is a good solution of the equation 
$$ 
A\circ \theta_{\f O}=\theta_{\f O}\circ F. \eqno{\Box}$$ 
\end{itemize}

\et

Primitive solution of \eqref{i1} such that $X$ belongs to the series $X=z^{n}$,   $X=T_{n}$, or $Z_n$ can be described explicitly
(see \cite{lattes}, Section 5). Below we will need only the following corollary of this description   (see \cite{lattes}, Corollary 5.2,  Corollary 5.5, and Corollary 5.8).

\bc \l{som1} Let $A,X,B$ be a primitive solution of \eqref{i1} such that $X=z^{n}$, where $n\geq 2.$ Then $\deg A\geq n$, unless 
$B=cz^{\pm m},$ $A=c^nz^{\pm m}.$ Similarly, if $X=T_{n}$, where $n>2,$ then $\deg A\geq n+1$, unless $B=\pm T_{m},$ $A=(\pm 1)^n T_{m}.$ Finally, if $X=Z_{n}$, where $n>2,$ then $\deg A\geq n+1$, unless 
$B=\v z^{\pm m},$ where $\v^{2n}=1,$ and $A=\v^n T_{m}.$ 
\ec

\end{subsection}

\end{section}

\begin{section}{Good chains}
Define a  {\it chain} $\f C=\f C(s,d)$ of length $s\geq 2$ as a sequence of  $s-1$ equalities 
\be \l{asd} U_{i}\circ V_{i}=V_{i+1}\circ U_{i+1},\ \ \ 1\leq i \leq s-1,\ee 
where  $U_i,$ $V_i,$ $1\leq i \leq s,$ are rational functions of degree at least two. 
Clearly, any such a chain corresponds to 
a sequence of  $s$ elementary transformations 
\be \l{cd} F\rightarrow F_1\rightarrow \dots \rightarrow F_s,\ee 
where   
$$F=V_1\circ U_1, \ \ \  F_i= U_i\circ V_i, \ \ \ \ \ 1\leq i \leq s.$$
The function $F=V_1\circ U_1$ is called the basis of  $\f C$, and the common degree $d$ of the functions in \eqref{cd} is called the degree of  $\f C$.

Any chain \eqref{asd} gives rise to  
the following commutative diagram:

\be 
\begin{CD} \l{gpa}
 @.  \vdots @. \vdots  @. \vdots @. \vdots\\
\hdots\ \ \  @. \C\P^1 @>U_2>> \C\P^1 @>U_3>> \C\P^1  @>U_4>> \C\P^1\\
@. @VV {V_1} V@VV {V_2} V @VV {V_3}  V @VV {V_4} V\\ 
@. \C\P^1 @>U_1>> \f \C\P^1 @>U_2>> \f\C\P^1@>U_3>> \C\P^1 \\
@. @.@VV {V_1} V @VV {V_2} V @VV {V_3} V\\ 
@. @.\C\P^1 @>U_1 >> \C\P^1  @>U_2>> \C\P^1\\
@. @. @. @VV {V_1} V @VV {V_2} V\\ 
 @. @. @.\C\P^1  @>U_1>> \C\P^1\, .\
\end{CD}
\ee 
\vskip 0.3cm
\noindent 
In particular, setting $$ V_{i,j}=V_{i}\circ V_{i+1}\circ \dots \circ V_{j}, \ \ \ \ U_{i,j}=U_{j}\circ U_{j-1}\circ \dots \circ U_i, \ \ \ 1\leq i\leq j\leq s,$$ 
we see that for any $i,j_1,j_2,$ $1\leq i\leq j_1\leq s,$ $1\leq i\leq j_2\leq s,$ the 
equality  
\be \l{eind} U_{i,j_1}\circ V_{i,j_2}=V_{j_1+1,j_1+j_2-i+1}\circ U_{j_2+1,j_1+j_2-i+1}\ee 
holds. Furthermore, 
the following statement is true (see \cite{lattes}, Lemma 3.1). 

\bl \l{leming} Let $\f C$ be a chain given by  \eqref{asd} and \eqref{cd}. 
Then  
$$V_{1,s}\circ U_{1,s}=F^{\circ s}, \ \ \ \ \  \ U_{1,s}\circ V_{1,s}=F_s^{\circ s}.\eqno{\Box}$$ 
\el

Let $\f C=\f C(s,d)$ be a chain.  Define its {\it dual} chain $\widehat{\f C}$ by the formulas
$$\widehat U_i=V_{s+1-i}, \ \ \ \ \widehat V_i=U_{s+1-i}, \ \ \ \ 1\leq i \leq s.$$
For  a natural number $k$ such that $$l_k=[s/k]\geq 2$$ define $\f C_k=\f C(l_k,d^k)$ as a chain corresponding to the sequence of $l_k-1$ equalities 
$$U_{1,k}\circ V_{1,k}=V_{k+1,2k}\circ U_{k+1,2k},$$ 
$$U_{k+1,2k}\circ V_{k+1,2k}=V_{2k+1,3k}\circ U_{2k+1,3k},$$ 
\vskip -0.2cm
$$U_{2k+1,3k}\circ V_{2k+1,3k}=V_{3k+1,4k}\circ U_{3k+1,4k},$$ 
$$\dots\, $$
$$U_{(l_k-2)k+1,(l_k-1)k}\circ V_{(l_k-2)k+1,(l_k-1)k}=V_{(l_k-1)k+1,l_kk}\circ U_{(l_k-1)k+1,l_kk}.$$ 

\vskip 0.2cm

A chain $\f C=\f C(s,d)$ is called {\it good} if all solutions of \eqref{m} provided by equalities \eqref{asd}  
 are good. 
For a good chain $\f C$ set 
$$d_1=\deg U_1=\deg U_2= \dots =\deg U_s,$$ and 
$$d_2=\deg V_1=\deg V_2= \dots =\deg V_s.$$ These  numbers are well defined  by Lemma \ref{good2},  and obviously satisfy the 
equality $d_1d_2=d.$ For good chains we will use the notation  $\f C=\f C(s,d_1,d_2)$ instead of the notation  $\f C=\f C(s,d).$
Clearly, Theorem \ref{sum+} implies inductively  
the following statement. 

\bl \l{lem0} Let  $\f C=\f C(s,d_1,d_2)$ be a good chain. Then any solution of \eqref{m} of the form \eqref{eind}  is good.  \qed
\el

In this  section we prove one of the 
main results of the paper: the finiteness of any good chain whose basis is non-special. 
More precisely, we prove the following.

\bt \l{mmtt}  Let  $\f C$ be a good chain of length $s$ and degree $d$ with non-special basis.   
Then $s < 12\log_{2}d+11 $.
\et

Notice that good chains of length $\approx\log_2d$ with non-special bases $F$ of degree $d$  exist and are easy to construct (see \cite{semi}, p. 1241). On the other hand, for special $F$  the theorem is not true. Indeed,  
 taking any commuting pair  $A$, 
$B$ of   powers, Chebyshev polynomials, or Latt\`es maps such that $\C(A,B)=\C(z)$, and setting $U_i = A,$ $V_i = B,$ $i \geq  1,$ we obviously obtain an infinite good 
chain.

\vskip 0.2cm Before proving Theorem \ref{mmtt} we will prove the following two lemmas.

\bl \l{som2} Let $A,C,D,B$ be a good solution of \eqref{m} such that either 
$C\underset{\mu}{\sim} T_n$ and $D\underset{\mu}{\sim} Z_{n/2}$, or  
$C\underset{\mu}{\sim} Z_{n/2}$ and $D\underset{\mu}{\sim}  z^{n}.$  Then $n\leq 2$. 
\el 
\pr 
Assume that $C\underset{\mu}{\sim} T_n$ and $D\underset{\mu}{\sim} Z_{n/2}$. 
 By Theorem \ref{t1}, $B:\f O_1^C\rightarrow \f O_1^D$ is a minimal holomorphic  map between orbifolds. 
On the other hand, if $n>2,$ then $\nu(\f O_1^{C})=\{2,2\}$, while the orbifold $\nu(\f O_1^{D})$ is non-ramified.
Since for such $\nu(\f O_1^{C})$ and $\nu(\f O_1^{D})$
  condition obviously 
 \eqref{rys} is not satisfied  at points $z$ where $\nu_{\f O_1^C}(z)=2,$ we conclude that $n\leq 2$ (in which case 
$\f O_1^{C}$ is non-ramified).

Assume now that $C\underset{\mu}{\sim} Z_{n/2}$ and $D\underset{\mu}{\sim}  z^{n}$.   
In this case $C$ and $D$ are the universal coverings of the orbifolds $\f O_2^C$ and $\f O_2^D$, so that
 $\Gamma_{\f O_2^C}=D_{n}$ and $\Gamma_{\f O_2^D}=C_{n}$.
Since $A:\f O_2^C\rightarrow \f O_2^D$ is a minimal holomorphic  map between orbifolds by Theorem \ref{t1}, it follows from Proposition \ref{poiu} that there exists a 
homomorphism $\phi:\, \Gamma_{\f O_2^C}\rightarrow \Gamma_{\f O_2^D}$ such that
$$ B\circ\sigma=\phi(\sigma)\circ B, \ \ \ \ \ \ \sigma\in \Gamma_{\f O_2^C}. $$  Moreover, if $n>2$, then $\Ker(\phi)\neq e,$ 
since $\vert D_{n}\vert =\vert C_{n}\vert $ but  the groups $D_{n}$ and $C_{n}$ are not isomorphic.
On the other hand,  $\Ker(\phi)$, as any other subgroup of  $\Gamma_{\f O_2^C}$, has the form $\Gamma_{\f O^{\prime}}$ for some orbifold $\f O^{\prime}$. Clearly,  the 
both functions $C$ and $B$ are invariant with respect to $\Gamma_{\f O^{\prime}}$, implying that they  
are rational functions in $\theta_{\f  O^{\prime}}$. Therefore, since $\Gamma_{\f O^{\prime}}\neq e$ implies that $\deg \theta_{\f  O^{\prime}}>1$, the quadruple 
$A,C,D,B$  is not a good solution of \eqref{m}. 
This contradiction shows that
$n\leq 2.$   \qed

We recall that the sets $E_1$ and $E_2$ are defined in Lemma \ref{xxx} and Lemma \ref{xxxx}.

\bl \l{xori}  Let $A$ be a rational function such that 
 $\chi(\f O_2^A)=0$, and $U,V$ rational functions of degree at least two such that 
 $A=U\circ V$ and $\deg U,\deg V\not\in E_2.$ Then \be \l{lama} \f O_1^A=\f O_1^V, \ \ \ \ \f O_2^V=\f O_1^U, \ \ \ \f O_2^A=\f O_2^U,\ \ \ \chi(\f O_2^V)=0.\ee

\el
\pr  Since $A:\f O_1^{A}\rightarrow \f O_2^{A}$ is a covering map between orbifolds, the maps
\be \l{azza} V:\f O_1^{A}\rightarrow U^*\,\f O_2^{A}, \ \ \ \ U:U^*\,\f O_2^{A}\rightarrow \f O_2^{A}\ee are covering maps  by Corollary \ref{indu2}. In particular, applying  \eqref{rhor} to the second map in \eqref{azza}, we conclude that 
 $$\chi(U^*\,\f O_2^{A})=0.$$ 
It follows now from Lemma \ref{xxxx} applied to covering maps \eqref{azza} that 
$$\f O_1^{A}=\f O_1^{V}, \ \ \ \ U^*\,\f O_2^{A}=\f O_2^{V},$$ 
and 
$$U^*\,\f O_2^{A}=\f O_1^{U}, \ \ \ \ \f O_2^{A}=\f O_2^{U},$$ implying that 
$$\f O_2^{V}=\f O_1^{U}. $$ 
Finally,  since $A:\f O_1^{A}\rightarrow \f O_2^{A}$ is a covering map, $\chi(\f O_1^A)=0$ by \eqref{rhor}. Thus, the first equality in \eqref{lama} implies that  $\chi(\f O_1^V)=0$, and 
applying \eqref{rhor} again we see that $\chi(\f O_2^V)=0$. \qed

\vskip 0.2cm
\noindent{\it Proof of Theorem \ref{mmtt}.}
Assume first that $d_1\geq d_2$ and, in addition, that 
$d_1$ is not contained in the set
$E_2$ from 
Lemma \ref{xxxx}. 
We show that 
in this case  \be \l{s} s<  \log_{d_2}84(d_1^3-2)+3 ,\ee unless the basis $F$ of $\f C$ is a special function.

Let us consider the equality 
\be \l{kr} U_{1,3}\circ V_{1,s-3} =V_{4,s}\circ U_{s-2,s}.\ee
Since
$$\deg U_{1,3}=d_1^3, \ \ \ \ \deg V_{1,s-3} =d_2^{s-3},$$
if
\be \l{inp} s\geq \log_{d_2}84(d_1^3-2)+3,\ee then 
$$\deg V_{1,s-3}=d_2^{s-3}\geq 84(d_1^3-2)= 84(\deg U_{1,3}-2).$$
Thus, since the solution of \eqref{m} provided by equality \eqref{kr}  is good by 
Lemma \ref{lem0}, it follows from Theorem \ref{m2} that $\chi(\f O_2^{U_{1,3}})\geq 0.$

Assume first that $\chi(\f O_2^{U_{1,3}})> 0.$ Since $d_1\not\in E_2$ implies that $d_1\neq 2$, and $d_1\neq 2$ implies that $d_1^3$ is 
not contained in the set
$E_1$,  Lemma \ref{xxx} implies  
that  $U_{1,3}$ is $\mu$-equivalent either to $z^{d_1^3}$, or $T_{d_1^3},$ or $Z_{d_1^3/2}$.
In case if $U_{1,3}$ is $\mu$-equivalent  to $z^{d_1^3}$, Lemma \ref{any} applied to the decompositions
$$U_{1,3}=U_{2,3}\circ U_1, \ \ \ \ \  U_{1,3}=U_{3}\circ U_{1,2}$$ 
implies that 
\be \l{llkk} U_{2,3}=\mu_2\circ z^{d_1^2}\circ \nu_1^{-1},\ \ \ U_1=\nu_1\circ z^{d_1} \circ \mu_1,\ee
$$ U_{3}=\mu_2\circ z^{d_1}\circ \nu_2^{-1},\ \ \ U_{1,2}=\nu_2\circ z^{d_1^2} \circ \mu_1,$$
for some M\"obius transformations $\mu_1,$ $\mu_2$, $\nu_1$, $\nu_2.$ 
Clearly, $d_1\geq d_2$ implies that 
$d_1^2> d_2$. Furthermore, the solution of \eqref{m} provided by the equality 
 \be \l{gops} U_{1,2}\circ V_1=V_3\circ U_{2,3}\ee is a good   by 
Lemma \ref{lem0}. 
Therefore, since
$$ U_{1,2}\underset{\mu}{\sim}z^{d_1^2},\ \ \ \ \ U_{2,3}\underset{\mu}{\sim} z^{d_1^2},  $$  we can apply 
Corollary \ref{som1} to \eqref{gops}, concluding that   
$$V_1=\mu_1^{-1}\circ cz^{\pm d_2}\circ \nu_1^{-1},  \ \ \ V_3=\nu_2\circ c^{d_1^2}z^{\pm d_2}\circ \mu_2^{-1},$$
for some $c\in \C$. Thus, 
$$F_1=U_1\circ V_1=\nu_1 \circ  c^{d_1}z^{\pm d_1d_2}\circ \nu_1^{-1}$$
is conjugate to $z^{\pm d_1d_2}.$  Since $F_1\sim F$, this implies by Lemma \ref{korova} that $F$ is special.
Similarly, if $U_{1,3}$ is $\mu$-equivalent  to $T_{d_1^3}$
we conclude that $F$ is conjugate to $\pm T_{d_1d_2}.$

Finally, the assumption that $U_{1,3}$ is $\mu$-equivalent  to $Z_{d_1^3/2}$ leads to a contradiction. Indeed, in this case 
applying Lemma \ref{any} inductively to the decomposition 
$$U_{1,3}=U_3\circ U_2\circ U_1,$$ we conclude that either 
$$U_3\underset{\mu}{\sim} Z_{d_1/2}, \ \ \  U_2\underset{\mu}{\sim} z^{d_1}, \ \ \  U_1\underset{\mu}{\sim} z^{d_1},$$ or
$$U_3\underset{\mu}{\sim} T_{d_1}, \ \ \  U_2\underset{\mu}{\sim} T_{d_1}, \ \ \  U_1\underset{\mu}{\sim} Z_{d_1/2},$$ or
$$U_3\underset{\mu}{\sim} T_{d_1}, \ \ \  U_2\underset{\mu}{\sim} Z_{d_1/2}, \ \ \  U_1\underset{\mu}{\sim} z^{d_1}.$$ 
Since the solution of \eqref{m} provided by the equality 
\be \l{ewa} V_3\circ U_{3}=U_{2}\circ V_2\ee is good, in the first case applying the second part of Lemma  \ref{som2} to \eqref{ewa} 
we obtain a contradiction with $d_1\neq  2$.
Similarly,  in the second case, we obtain a contradiction applying the first part of Lemma  \ref{som2}  to the equality 
\be \l{ewa0} V_2\circ U_{2}=U_ {1}\circ V_1.\ee  
Lastly, in the third case we obtain a contradiction applying Lemma   \ref{som2} to either of equalities \eqref{ewa}, \eqref{ewa0}.

Assume now that $\chi(\f O_2^{U_{1,3}})= 0.$ Since $d_1\not\in E_2$ implies $d_1^2\not\in E_2$, it follows from 
Lemma \ref{xori} applied to the decomposition $U_{1,3}=U_3\circ U_{1,2}$
that \be \l{0x} \chi(\f O_2^{U_{1,2}})= 0.\ee Applying now Lemma \ref{xori} to the decomposition   $U_{1,2}=U_2\circ U_1$, 
we obtain the equalities 
\be \l{jo2}\f O_2^{U_ {1}}=\f O_1^{U_{2}},\ee 
\be \l{0y} \f O_2^{U_{1,2}}=\f O_2^{U_2},\ee
\be \l{2jo}\chi(\f O_2^{U_ {1}})=0.\ee 
Clearly, \eqref{0x} and \eqref{0y} imply that 
\be \l{0z} \chi(\f O_2^{U_2})=0.\ee 
Further, since Theorem \ref{t1} applied to \eqref{ewa0} implies that 
$V_2: \f O_2^{U_2}\rightarrow \f O_2^{U_1}$ is a minimal holomorphic map between orbifolds, it follows from  equalities  \eqref{2jo} and \eqref{0z}
by Proposition \ref{p1} that $V_2: \f O_2^{U_2}\rightarrow \f O_2^{U_1}$ is a covering map between orbifolds.
Since $U_2: \f O_1^{U_2}\rightarrow \f O_2^{U_2}$ also  is a covering map,
it follows now from Corollary \ref{serka0} that 
$$F_1=V_2\circ U_2:\f O_1^{U_2}\rightarrow \f O_2^{U_1}$$ 
 is a covering map too. Therefore,  in view of equality 
 \eqref{jo2}, the function $F_1$ is a Latt\`es map, implying by Lemma \ref{korova} that the function $F$ is special.

We proved that under the assumptions $d_1\geq d_2$ and
$d_1\not \in E_2$ equality \eqref{s} holds, unless $F$ is a special function. 
Let us explain now how to get rid of these assumptions and obtain the inequality from the formulation of the theorem. 
First, if $d_1$ is contained in the set
$E_2$, we can consider instead of the chain $\f C=\f C(s,d_1,d_2)$ the chain 
$\f C_4=\f C(\left[\frac{s}{4}\right],d_1^4,d_2^4)$, which is also good by Lemma \ref{lem0}.  Since for any number $d_1\geq 2$ the number $d_1^4$ does not belong to $E_2$, the above argument shows that 
\be \l{pqa} \left[\frac{s}{4}\right]<  \log_{d_2^4}84(d_1^{12}-2)+3,\ee
unless the basis $\t F$ of $\f C_4$ is a special function. On the other hand, since $$\t F=V_{1,4}\circ U_{1,4}=F^{\circ 4}$$ by Lemma \ref{leming}, it follows 
from Lemma  \ref{baran}  that if  $\t F$ is a special function, then $F$  also is  a special function.

Furthermore, since 
$$d_2^4\geq 2^4, \ \ \ \ d_1^{12}\leq \left(\frac{d}{2}\right)^{12},$$ and \eqref{pqa} implies that
$$\frac{s}{4}<  \log_{d_2^4}84(d_1^{12}-2)+4,$$
we conclude that inequality \eqref{pqa} yields the inequality
$$s< 4\log_{d_2^4}84(d_1^{12}-2)+16\leq \log_{2}84\left(\left(\frac{d}{2}\right)^{12}-2\right)+16< $$
$$< \log_{2}84+\log_{2}\left(\frac{d}{2}\right)^{12}+16< 12\log_{2}d+11$$
from the formulation of the theorem. Thus, if the inequality \be \l{pes} s< 12\log_{2}d+11\ee
does not hold, inequality \eqref{pqa} does not hold either, implying that $F$ is special.

Finally,  if $d_1< d_2$,  we can consider instead of the chain $\f C=\f C(s,d_1,d_2)$ its dual chain $\hat{\f C}=\hat{\f C}(s,d_2,d_1)$. By the above argument, inequality \eqref{pes} holds, 
unless  the basis $\hat F$ of $\hat{\f C}$ is a special function. 
On the other hand, since $$\hat F=U_s\circ V_s\sim F,$$ if the basis $\hat F$ of $\hat{\f C}$ is a special function, then the basis $F$ of $\f C$ is  also  a special function by 
Lemma \ref{korova}. 
\qed 
\end{section}

\begin{section}{Group $G(B)$ and primitive solutions} 
For a rational function $B$, we denote by $\f E(B)$  the set of rational functions $X$ of degree at least two such that \eqref{i1} holds for some rational function $A$, and by $\f E_0(B)$ the subset of $\f E(B)$ consisting of functions $X$
such that $\C(X,B)=\C(z)$. In addition, we denote by $\overline{\f E_0(B)} $ the quotient set of $\f E_0(B)$ by the equivalence relation 
which identifies $X_1,$ $X_2\in \f  E_0(B)$ if there exists a M\"obius transformation $\mu$ such that 
\be \l{meb} X_1=\mu \circ X_2.\ee 
In this section, we prove the finiteness of the set $\overline{\f E_0(B)} $ 
for non-special $B$. Abusing the notation we will denote by $X$ both an element of $\f E_0(B)$ and its equivalence class in 
$\overline{\f E_0(B)} $.


We recall that the symmetry group of a rational function $F$ is defined as the group of all M\"obius transformations $\mu$ commuting with $F$. Since such transformations map 
periodic points of $F$ of any given period to themselves and any  M\"obius transformation is defined by 
its values at any three points, the symmetry group of any rational function is finite.   
We define {\it the extended symmetry group }
$G(F)$ of $F$ as the group of M\"obius transformations $\sigma$ such that \be \l{eblys} F\circ \sigma=\nu_{\sigma} \circ F\ee for some  M\"obius transformations $\nu_{\sigma}$.  It is easy to see that 
 $G(F)$ is indeed a group with respect to the composition operation and that the map $$\gamma_F:\sigma \rightarrow \nu_{\sigma}$$ is a homomorphism from   $G(F)$ to the group 
$\rm{Aut}(\C\P^1)$. 

We denote by $\f D$ the subgroup of $\rm{Aut}(\C\P^1)$  consisting of the transformations $\sigma= c z^{\pm 1}$, $c\in \C\setminus\{0\}$. Notice that $\f D$ can be described as the subgroup of $\rm{Aut}(\C\P^1)$  consisting 
of all M\"obius transformations $\sigma$ such that $\sigma\{0,\infty\}=\{0,\infty\}$, or 
equivalently such that $\sigma^{-1}\{0,\infty\}=\{0,\infty\}$.

\bl \l{zn} For $F=z^{\pm d}$ the group $G(F)$ is  $\f D.$
\el 
\pr It is clear that $\f D\subseteq G(F)$. On the other hand, 
if equality \eqref{eblys} holds, then  $\nu_\sigma^{-1}\{0,\infty\}=\{0,\infty\}$, since otherwise the preimage 
$(\nu_{\sigma}\circ  F)^{-1}\{0,\infty\}$ and hence 
 the preimage $(F\circ \sigma)^{-1}\{0,\infty\}$ contains more than two points.
Therefore,  
for any $\sigma\in G(F)$ the transformation $\nu_{\sigma}$ belongs to $ \f D$. Now \eqref{eblys} implies that $\sigma$ also belongs to $\f D$.  
\qed

\vskip 0.2cm

\bt \l{prim} Let $F$ be a rational function of degree $d\geq 2$ such that  $F\underset{\mu}{\not\sim} z^d.$  Then the group $G(F)$ is isomorphic to one of the five finite rotation groups of the sphere $A_4,$ $S_4,$ $A_5,$ $C_n$, $D_{2n}$,
and	the order  of any element of $G(F)$ does not exceed  $d.$ In particular, $\vert G(F)\vert \leq \max\{60,2d\}.$
\et 

\pr 
Any non-identical element of the group $\Aut(\C\P^1)\cong \rm{PSL}_2(\C)$ is conjugate either to $z\rightarrow z+1$ or to $z\rightarrow \lambda z$ for some $\lambda \in \C\setminus\{0,1\}.$ Thus, making the change 
\be \l{zve} F\rightarrow \mu_1\circ F\circ \mu_2, \ \ \ \sigma\rightarrow \mu_2^{-1}\circ \sigma\circ \mu_2,\ \ \ \nu_{\sigma}\rightarrow \mu_1\circ \nu_{\sigma}\circ \mu_1^{-1}\ee for convenient $\mu_1,$ $\mu_2\in\rm{Aut}(\C\P^1)$, without loss of generality we may assume that $\sigma$ and $\nu_{\sigma}$ in \eqref{eblys} have one of the two forms above.  

We observe first that 
the equalities  \be \l{j1} F(z+1)=\lambda F(z), \ \ \ \lambda\in \C\setminus\{0,1\},\ee and  \be \l{j2}F(z+1)=F(z)+1\ee 
are impossible. Indeed, if $F$ has a finite pole, then any of these equalities implies that $F$ 
has infinitely many poles. On the other hand, if $F$ is a polynomial of degree $d\geq 2$, then we obtain a contradiction comparing the coefficients of $z^d$ in the left and the right sides of 
\eqref{j1}, and the coefficients of
 $z^{d-1}$ in left and the right sides of   \eqref{j2}, correspondingly. 

Furthermore, comparing the free terms in  the Laurent series at infinity of the left and the right sides of the equality \be \l{egga} F(\lambda z)=F(z)+1, \ \ \ \lambda\in \C\setminus\{0,1\},\ee 
we conclude that this equality is impossible either. 
Thus, \be \l{zore} F(\lambda_1z)=\lambda_2F(z),\ \ \ \lambda_1, \lambda_2\in \C\setminus\{0,1\}.\ee Comparing now coefficients in the left and the right sides of \eqref{zore} and taking into account that 
$F\underset{\mu}{\not\sim} z^d$, we conclude that $\lambda_1$ is a root of unity.
Furthermore, the order of the transformation $z\rightarrow \lambda_1z$ in the group  $G(F)$ does not exceed the maximum number $n$ such that $F$ can be represented in the form \be \l{repr} F=z^rR(z^n), \ \ \ \ R\in \C(z).\ee In particular, the order of any 
element of $G(F)$ does not exceed $d.$
 Indeed, since $F\underset{\mu}{\not\sim} z^d$, the function $R$ in \eqref{repr} has a zero or a pole  distinct from $ 0$ and $\infty$, implying that   $d\geq n$.

The above argument shows that any element of $G(F)$ has finite order. In order to prove 
the finiteness of $G(F)$ we will use  the Schur theorem,  which states that if any element of a finitely generated subgroup $G$ of  $\rm{GL}_k(\C)$ has finite order, then $G$ has finite order (see e.g. \cite{cur}, (36.2)).
Specifically, assume that $G(F)$ is infinite, and let  $\sigma_1,\sigma_2, \dots ,\sigma_s, \dots$ be an infinite sequence of pairwise distinct elements of $G(F)$. 
Observe first that for any $s\geq 1$ the finitely generated group  $$\Gamma_s=<\sigma_1,\sigma_2, \dots, \sigma_s>$$ is finite.  
Indeed, if $\Gamma_s$ is infinite, then its lifting $$\widehat{\Gamma}_s\subset \rm{SL}_2(\C)\subset  \rm{GL}_2(\C)$$ is  also 
infinite, implying by the Schur theorem that $\widehat{\Gamma}_s$ has 
an element of  infinite order. But  in this case $\Gamma_s$ also  has 
an element of  infinite order  in contradiction with the fact that any element  of $G(F)$ has finite order.

Since the elements  $\sigma_1,\sigma_2, \dots ,\sigma_s, \dots$  are pairwise distinct, 
$\vert \Gamma_s\vert \to \infty$. On the other hand, 
since the groups $\Gamma_s,$ $s\geq 1,$ 
are finite subgroups of $\Aut(\C\P^1)$, they belong to the list $A_4,$ $S_4,$ $A_5,$ $C_n$, $D_{2n}$. Therefore,  for $s$ big enough the group $\Gamma_s$ is either $C_n$ or $D_{2n}$, where $n$ satisfies the inequality $n>d$. However, since the both groups $C_n$ and $D_{2n}$ have an element of order $n$, this  contradicts to the fact that the order of any element  of $G(F)$ does not exceed $d.$
Thus, the group $G(F)$ is finite. 
Finally, if $G(F)$ is 
$A_4,$ $S_4,$ or $A_5,$ then $G(F)\leq 60$, while if $G(F)$ is 
$C_n$ or $D_{2n}$, then $n\leq d$, since $C_n$ and $D_{2n}$ have an element of order $n.$   
\qed

\bl \l{meds} Let $\delta$ be a M\"obius transformation that does not belong to  $\f D$.
Then the intersection of 
the groups $\f D$ and $\delta^{-1}\circ \f D\circ \delta$ 
is a finite group isomorphic to a  subgroup of
 the Klein four-group. 
\el 
\pr Let us denote the intersection of 
the groups $\f D$ and $\delta^{-1}\circ \f D\circ \delta$ by  $\f R$. 
We show first that $\f R$ consists of involutions.  
Assume in contrary that  $\mu\in \f R$ is not an involution, and let $\mu'$ be an element of $\f D$ which makes the diagram 
\be \l{dd1} 
\begin{CD}
\C\P^1 @>\mu>> \C\P^1\\
@VV\delta V @VV\delta V\\ 
 \C\P^1 @>\mu' >> \C\P^1\ 
\end{CD}
\ee
commutative. Since the diagram 
\be \l{dd2}
\begin{CD}
\C\P^1 @>\mu^{\circ 2}>> \C\P^1\\
@VV\delta V @VV\delta V\\ 
 \C\P^1 @>\mu'^{\circ 2} >> \C\P^1\ 
\end{CD}
\ee 
also is commutative, $\mu'$ is not an involution either. Therefore, since  
$$\mu\{0,\infty\}=\{0,\infty\}, \  \ \ \ \mu'\{0,\infty\}=\{0,\infty\},$$  the set $\{0,\infty\}$ is the fixed point set of the transformations $\mu^{\circ 2}$ and $\mu'^{\circ 2}$. It follows now from \eqref{dd2} that 
$\delta\{0,\infty\}=\{0,\infty\}$, in contradiction with $\delta\not\in \f D.$

Since $\f R$ consists of involutions, any $\mu\in \f R$ has either the form $\mu=\pm z,$  or the form $\mu=c/z,$ $c\in \C\setminus\{0\}$. Furthermore, if  two transformations  $\mu_1=c_1/z,$ $c_1\in \C\setminus\{0\}$, and $\mu_2=c_2/z,$ $c_2\in \C\setminus\{0\}$, belong to $\f R$, then  their composition $\frac{c_1}{c_2}z$ also belongs to $\f R$ implying that 
$c_2=\pm c_1.$ Therefore, the group $\f R$ contains at most four elements: $\pm z,$ $\pm \frac{c_1}{z}$.  \qed

\bt \l{karas}
Let $B$ be a  non-special rational function of degree $d\geq 2$. 
Then for any $X\in\f E_0(B)$ the inequalities $\chi(\f O_2^{X})>0$ and $\deg X\leq \max\{60,2d\}$ hold. Furthermore, 
the set $\overline{\f E_0(B)} $ is finite and its cardinality can be bounded from above by a number depending on $d$ only.
\et

\pr Let $X$ be an element of $\f E_0(B)$ and $A$ the corresponding function such that \eqref{i1} holds.
By Theorem \ref{oip},  the diagram
\be  \l{lepss}
\begin{CD}
\f O_1^X @>B>> \f O_1^X\\
@VV X V @VV X V\\ 
\f O_2^X @>A >> \f O_2^X\ 
\end{CD}
\ee
consists of minimal holomorphic  maps between orbifolds, and  $\chi(\f O_2^{X})\geq 0$. 
Moreover, in fact $\chi(\f O_2^{X})> 0.$ Indeed, 
if $\chi(\f O_2^{X})= 0,$ then  $\chi(\f O_1^{X})= 0$ by \eqref{rhor}. Therefore, since 
\be \l{ebys} 
B\, : \f O_1^{X}\rightarrow \f O_1^{X}\ee is a minimal holomorphic map, it follows from Proposition \ref{p1} that  \eqref{ebys} is a covering map, implying that 
$B$ is a Latt\`es map, in contradiction with the assumption. Notice that $\chi(\f O_2^{X})> 0$ implies $\chi(\f O_1^{X})> 0$ by \eqref{rhor}.

Since $B$ is non-special and \eqref{ebys} is a minimal holomorphic map, it follows from Theorem \ref{uni}  
that the orbifold $\f O^B_0$ is defined and ${\f O}_1^X\preceq \f O^B_0$. Moreover, since $B$ is not a Latt\`es map, $\chi(\f O_0^{B})> 0.$ We observe first that the number of orbifolds $\f O$ such that $\f O\preceq \f O^B_0$ is finite and can be bounded by a number depending on $d$ only. Clearly, it is enough to show that if  
 $\f O^B_0$ belongs to the series  $\{n,n\}$, $n\geq 2,$ or $\{2,2,n\}$, $n> 2,$ 
then  $n$ is bounded in terms of $d.$ By Theorem \ref{las}, there exists a rational function $F_0$  such that  the  diagram 
\be \l{xcv} 
\begin{CD} 
\C\P^1 @>F_0>> \C\P^1 \\
@VV \theta_{\f{O}_0^B} V @VV \theta_{\f{O}_0^B} V\\ 
\f{O}_0^B @>B>> \f{O}_0^B \,.
\end{CD}
\ee commutes and the triple $B,$ $\theta_{\f O_0^B}$,   $F_0$  is a primitive solution of equation \eqref{i1}.
By condition, $B$ is not special.   Therefore, if $\theta_{\f{O}_0^B}$ is $\mu$-equivalent  to $z^n$ or to $Z_n$, applying Corollary \ref{som1} to \eqref{xcv}   we conclude that $n\leq d$.  

Since the number of orbifolds $\f O$ such that $\f O\preceq \f O^B_0$ is finite, in order to prove the finiteness of $\overline{\f E_0(B)}$ it is enough to show that  for any fixed orbifold $\f O$ with $\chi(\f O)>0$ there exist only finitely many $X \in \overline{\f E_0(B)}$ 
 such that $\f O_1^X=\f O.$ 
Assume first that $\f O$ is non-ramified. Then 
$X=\theta_{\f O_2^X}$  by \eqref{ravv}, and diagram \eqref{lepss} reduces to diagram \eqref{dia3}. Furthermore, by Theo\-rem \ref{las},  there exists an  automorphism
$\phi:\,\Gamma_{\f O_2^X}\rightarrow \Gamma_{\f O_2^X}$ such that 
for any $\sigma\in \Gamma_{\f O_2^X}$ the equality 
\be \l{hol} B\circ\sigma=\phi(\sigma)\circ B\ee holds. 
 Therefore, $\Gamma_{\f O_2^X}$ belongs to the intersection 
$$G_0(B)= G(B)\cap \gamma_B(G(B)),$$ and hence 
the number of classes $X$ in $\overline{\f E_0(B)}$ with non-ramified $\f O_1^X$
does not exceed the number of  subgroups of 
$G_0(B)$. Since $G_0(B)\subseteq G(B)$ and $$\vert G(B)\vert \leq \max\{60,2d\}$$ whenever  $B\underset{\mu}{\not\sim}z^d$, we conclude that for $B\underset{\mu}{\not\sim}z^d$ the number of subgroups of  $G_0(B)$ and therefore  the cardinality of $\overline{\f E_0(B)}$ is bounded from above by a number depending on $d$ only. Moreover, 
\be \l{rad} \deg X=\deg \theta_{\f O_2^X}=\vert\Gamma_{\f O_2^X}\vert \leq \vert G(B)\vert\leq \max\{60,2d\}.\ee

Assume now that $B\underset{\mu}{\sim}z^d$. 
Without loss of generality we may assume that $B=\delta\circ z^d,$ where $\delta$ is a M\"obius transformation, so that $G(B)=\f D$ by Lemma \ref{zn}, and $$\gamma_B(G(B))=\delta\circ \f D\circ \delta^{-1}.$$ 
By Lemma \ref{meds}, either $G_0(B)$ is a group of order at most four, or 
$\delta\in \f D$. 
In the first case, the number of classes $X$ in $\overline{\f E_0(B)}$ with non-ramified $\f O_1^X$ does not exceed four,
 and  inequality \eqref{rad} still holds. 
On the other hand, in the second case $B=cz^{\pm d}$, $c\in \C\setminus\{0\},$ implying that $B$ is conjugate to $z^{\pm d}$ in 
contradiction with the assumption.

The case where $\f O=\f O_1^X$ is ramified reduces to the previous one. Indeed, since \eqref{ebys} is a minimal holomorphic map, it follows from Theorem \ref{las} that there exists a rational function $F$  
such that  for any $X$ with $\f O_1^X=\f O$ 
the  diagram 
\be 
\begin{CD} \l{gooopa3}
\C\P^1 @>F>> \C\P^1 \\
@VV \theta_{\f{O}_1^X} V @VV \theta_{\f{O}_1^X} V\\ 
\f{O}_1^X @>B>> \f{O}_1^X\\
@VV X V @VV X V\\ 
{\f O}_2^X @>A >> \ \ {\f O}_2^X\,,\  
\end{CD}
\ee
commutes, and the triple $B,$ $\theta_{\f O_1^X}$,   $F$  is a primitive solution of equation \eqref{i1}. Set $$X'=X\circ \theta_{\f O_1^X}.$$ Since 
$X'=\theta_{\f O_2^X}$ by \eqref{ravv}, it follows from the commutativity of diagram \eqref{gooopa3} by Theorem \ref{sum+} 
that any $X\in \overline{\f E_0(B)} $ with $\f O_1^X=\f O$ lifts to some $X'\in \overline{\f E_0(F)}$ with non-ramified
$\f O_1^{X'}.$
Moreover, since $B$ is not a special function, $F$ also is not a special function, 
 by Lemma \ref{zaq}.   
Since $\deg F=\deg B,$ this yields that  the number of classes $X$ in $\overline{\f E_0(B)}$ with $\f O_1^X=\f O$ can be bounded from above by a number depending on $d$ only.
Lastly, the inequality  $$\deg X\leq \max\{60,2d\}$$ still holds since
$$\deg X<\deg X'=\deg \theta_{\f O_2^X}.\eqno{\Box}$$
\end{section}

\begin{section}{Description of $\f E(B)$ for non-special $B$}
In this section, we  describe the structure of the set $\f E(B)$ for  non-special rational functions $B$, and prove  Theorem \ref{main1} and Theorem \ref{main2}. 
We also prove an effective version of the Ritt theorem.

\bt \l{eog}
Let $B$ be a  non-special rational function of degree $d\geq 2$.
Then any $X\in \f E(B)$ can be represented in the form \be \l{repres} X=X'\circ U\circ B^{\circ k},\ \ \ \ k\geq 0,\ee where $\chi(\f O_2^{X'})> 0$ and $U$ is a compositional right factor of some iterate $B^{\circ l}$, $l\geq 0.$
Furthermore,  \be \l{nerr} \deg X'\leq \max\{60,2d\}\ee and 
\be \l{pipez} l\leq  (\sigma_0(d)-2)(12\log_{2}d+11),\ee 
where $\sigma_0(d)$ is the number of divisors function. 
\et 
\pr If $\C(X,B)=\C(z)$, then the statements of the theorem holds for $l=0$ and $k=0$ by Theorem \ref{karas}. So, assume that 
$\C(X,B)\neq\C(x)$. Without loss of generality we may assume that $X$ is not a rational function in $B^{\circ k}$, $k\geq 1,$ since if $X=\widehat X\circ B^{\circ k}$ is contained in  $\f E(B)$, then $\widehat  X$ also is contained in $\f E(B).$

Recall that any solution $A,X,B$ of \eqref{i1} can be reduced to a primitive one as follows.  
Let $\C(X,B)=\C(U_1)$, where $U_1\in \C(z)$,  and let $X_1,V_1$ be rational functions such that 
\be \l{off} X=X_1\circ U_1,\ \ \  B=V_1\circ U_1,  \ee 
and $\C(X_1,V_1)=\C(z).$ Substituting expressions \eqref{off} in \eqref{i1} we see  that the diagram 
$$
\begin{CD} 
\C\P^1 @> B=V_1\circ U_1>>\C\P^1 \\ 
@V U_1  VV @VV U_1   V\\ 
 \C\P^1 @>  U_1\circ V_1>> \C\P^1
 \\ 
@V {X_1}  VV @VV {X_1}   V\\ 
 \C\P^1 @> A>> \C\P^1, 
\end{CD} 
$$
commutes. Thus, 
$$A'=A, \ \ \ \ X'=X_1, \ \ \ \ B'=U_1\circ V_1$$ is also a solution  of \e{i1}.
The new solution   is not necessary primitive. Nevertheless, $\deg X_1 < \deg  X.$ Therefore,
after a finite number of similar transformations we will arrive to a primitive
solution.  

In more details, we can find 
rational functions $X_{i},V_{i},$  $U_{i}$, $1\leq i\leq l$, such that
\be \l{tog1} \C(X_i,U_i\circ V_i)=\C(U_{i+1}),  \ \ \ \ \deg U_{i+1} >1, \ \ \ 1\leq i \leq l-1,\ee 
\be \l{tog2} X_i=X_{i+1}\circ U_{i+1},\ \ \ \ U_i\circ V_i=V_{i+1}\circ U_{i+1}, \ \ \ 1\leq i \leq l-1,\ee
\be \l{ex} \C(X_{i},V_{i})=\C(z),  \ \ \ 1\leq i \leq l,\ee 
and 
\be \l{tog4} \C(X_l,U_l\circ V_l)=\C(z).\ee
Setting 
\be \l{udovv} B_i= U_i\circ V_i=V_{i+1}\circ U_{i+1},\ \ \ 1\leq i \leq l-1,\ee
and   $$U=U_l\circ \dots \circ U_2\circ U_1, $$
we see that by construction 
$$X=X_l\circ U,$$ 
 the diagram 
\be \l{epta}
\begin{CD} 
\C\P^1 @> B>>\C\P^1 \\ 
@V U  VV @VV U   V\\ 
 \C\P^1 @> B_l>> \C\P^1
 \\ 
@V X_l  VV @VV X_l   V\\ 
 \C\P^1 @> A>> \C\P^1 
\end{CD} 
\ee
commutes, and 
$A,X_l, B_l$ is a primitive solution of \eqref{i1}. 
Thus, since $U$
is a compositional right factor of $B^{\circ l}$ by Lemma \ref{leming}, and the inequalities $\deg X_l\leq \max\{60,2d\}$ and $\chi(\f O_2^{X_l})>0$ hold by Theorem \ref{karas},
 in order to prove the theorem we only must show that  $l$ 
satisfies  inequality \eqref{pipez}.

Since the first equality in \eqref{tog2} implies that any  compositional right factor of $U_{i+1}$ is a
compositional right factor of $X_i$, it follows from \eqref{ex} that 
\be \l{vot} \C(V_i,U_{i+1})=\C(z), \ \ \ 1\leq i \leq l-1.\ee
In  turn, \eqref{vot} yields that 
\be \l{ner} \deg U_{i+1}\leq \deg U_i, \ \ \ i\geq 1.
\ee
Indeed, 
denote by $ I$ the imprimitivity system  of the monodromy group of $B_i$ corresponding to the decomposition 
$B_i=V_{i+1}\circ U_{i+1}$, and by $J$ the imprimitivity system   corresponding to the decomposition 
$B_i= U_i\circ V_i$. Since  each block of  $ I$  contains $\deg U_{i+1}$ elements, while the number of blocks 
of  $J$ is equal to $\deg U_{i}$, if $\deg U_{i+1}> \deg U_i$, then 
there is a block of $J$ containing at least two elements from a block of $I,$
implying that $\C(V_i,U_{i+1})\neq \C(z).$

By assumption, $\deg U_i\geq 2$, $1\leq i \leq l.$ On the other hand, since $X$ is not a rational function in $B$ the inequality $\deg U_1<d$ holds,  implying by \eqref{ner} that $\deg U_i< d$, $1\leq i \leq l.$ Thus, the functions $V_{i},$  $U_{i}$, $1\leq i\leq l$,
have degree at least two and hence the sequence 
\be \l{tog3}  U_i\circ V_i=V_{i+1}\circ U_{i+1}, \ \ \ 1\leq i\leq l-1, \ee 
is a chain. We denote this chain by $\f C$. 
Let 
$$1\leq k_1\leq \dots \leq k_r\leq  l$$
 be  indices  such that
\be \l{indi}  \deg U_1=\deg U_2=\dots =\deg U_{k_1}, \ee
$$\deg U_{k_1+1}=\deg U_{k_1+2}=\dots =\deg U_{k_2}, $$
\vskip -0.6cm
$$\dots$$
$$ \deg U_{k_{r}+1}=\deg U_{k_{r}+2}=\dots =\deg U_{l}, $$
and
\be \l{steve}\deg U_ {1}>\deg U_{k_1+1} >\deg U_{k_2+1}> \dots >  \deg U_{k_{r}+1} .\ee 
Setting for convenience 
$k_0=0$, $k_{r+1}=l$, we see that in view of conditions \eqref{vot} and \eqref{indi},  for any $j,$ $0\leq j \leq r,$ the subchain 
\be \l{sl} U_i\circ V_i=V_{i+1}\circ U_{i+1}, \ \ \ k_j+1\leq i\leq k_{j+1}-1,\ee
of the chain $\f C$ is good by Lemma \ref{good}. Therefore, by Theorem \ref{mmtt}, its length $k_{j+1}-k_j$ is less than or equal to $12\log_{2}d+11$, unless its basis 
is special. On the other hand, 
since by construction 
$$B\rightarrow B_1\rightarrow B_2\rightarrow \dots \rightarrow B_l$$
is a chain of elementary transformations, it follows from  
Lemma \ref{korova} that the  basis of  \eqref{sl}
is special if and only if $B$ is special. Thus, we conclude that 
$$l=\sum_{j=0}^{r}(k_{j+1}-k_j)\leq (r+1)(12\log_{2}d+11),$$ 
unless $B$
is special.
Finally, since  $\deg U_i$, $1\leq i \leq l,$ is a proper divisor of $d$, it follows from \eqref{steve} that
that $$r+1\leq \sigma_0(d)-2.\eqno{\Box}$$

\vskip 0.2cm

\vskip 0.2cm 
\noindent{\it Proof of Theorem \ref{main1}.} Assume that $X$ is an element of $\f E(B)$ such that 
$X$ is not a rational function in 
$B^{\circ k}$, $k\geq 1,$ and consider diagram \eqref{epta}. 
Since decompositions of a rational function $R$ of degree $n$ into a composition of rational functions $R=L\circ M$, considered up to the change 
\be \l{chna} L= L\circ \nu^{-1}, \ \ \ \ \ M=\nu\circ M,\ee
correspond to imprimitivity systems of the monodromy group ${\rm Mon}(R)\subset S_n$ of $R$,
 there 
exists a  function $\omega:\N\rightarrow \N$ such that, up to the change \eqref{chna}, any rational function  of degree $n$
has at most $\omega(n)$ decompositions into a composition of rational functions.  
Thus, since $U$
is a compositional right factor of $B^{\circ l}$, where $l$ satisfies  \eqref{pipez},  there exist $W_1,$ $W_2,\dots ,W_{s}$, where $s$ can be bounded from above by a number depending on $d$ only,
such that $U$ in \eqref{epta} has the form $U= \nu\circ W_j$ for some $j,$ $1\leq j \leq s$, and a M\"obius transformation $\nu$.
Besides, it is clear that 
$$\deg W_j \leq \deg B^{\circ l}=d^l, \ \ \ \ 1\leq j \leq s.$$

Changing in diagram \eqref{epta} the function $X_l$ to the function $X_l\circ \nu^{-1}$, the function $U$ to the function $\nu\circ U$, and the function $B_l$ to the function $\nu\circ U\circ \nu^{-1}$, for a convenient M\"obius transformation $\nu$, without loss of generality we may assume that   $U=W_j$, for some $j,$ $1\leq j \leq s$. The function $B_l$ is defined in a unique way by $U$, and it follows from Theorem \ref{karas} that for each $B_l$,  up to the change $X_l\rightarrow \mu\circ X_l$, where $\mu$ is a M\"obius transformation,  there exist only finitely many functions $X_l$ satisfying \eqref{epta}. Moreover,  the number of such functions and their degrees  can 
be bounded from above in terms of $d$ only.
 Finally, it is clear that if to the function $X_l$ corresponds the function $A,$ then to the function $\mu\circ X_l$
corresponds the conjugate function $\mu\circ A\circ \mu^{-1}.$ \qed 

\vskip 0.2cm

\br 
Since $A\sim B$ implies that $B$ is semiconjugate to $A$, 
Theorem \ref{main1} yields in particular that if $B$ is not special, then $[B]$ contains only finitely many conjugacy classes. 
This result also follows from the McMullen theorem \cite{Mc} about isospectral rational functions  (see  \cite{rec}). However, the McMullen theorem is non-effective, while Theorem \ref{main1} asserts that  the number of  classes in $[B]$ is bounded from above by a number depending on $\deg B$ only.  

\er

\noindent{\it Proof of Theorem \ref{main2}.} Assume that $X$ commutes with $B.$ Without loss of generality we may assume 
that $X$ is not a rational function in $B^{\circ k}$, $k\geq 1,$ since if $$X=\widehat X\circ B^{\circ k}$$ commutes with $B$, then $\widehat  X$ also 
 commutes with $B$.
By Theorem \ref{main1}, there exist  $X_1, X_2,\dots, X_{r}\in \f E(B)$  
such 
that $X=\mu \circ X_j$ for some $j,$ $1\leq j \leq r,$ and a M\"obius transformation $\mu$. Therefore, in order to prove Theorem \ref{main2}, it is enough to show that for any $X\in \f E(B)$ 
there exist only finitely many $\mu$ such that $\mu\circ X$ commutes with $B.$ 

Take arbitrary $\mu_0$ such  $\mu_0\circ X $ commutes with $B$, and assume that  
$\mu\circ X $ also commutes with $B$.
We have:
$$B\circ \mu_0\circ X =\mu_0\circ X \circ B=(\mu_0\circ \mu^{-1})\circ \mu\circ X \circ B=
 (\mu_0\circ \mu^{-1})\circ B\circ \mu\circ X .$$ 
Therefore,
$$B\circ \mu_0= (\mu_0\circ \mu^{-1})\circ B\circ \mu.$$ 
implying that 
 $\nu= \mu_0 \circ \mu^{-1}$ commutes with $B$.  Since the number of M\"obius transformation  commuting with $B$ is finite, this implies the required statement. \qed

\vskip 0.2cm

In conclusion, we prove the following effective version of the Ritt theorem (cf. \cite{r}, \cite{e2}).

\bt \l{main3} Let  $B$ and $X$ be commuting rational functions of degree at least two. If $B$ is not special, then $X^{\circ l}= B^{\circ k}$ for some $l,k\geq1.$ Furthermore, there exists a (computable) function $\delta:\N\rightarrow \N$ such that for any non-special $B$ of degree $d$ the number $l$ is bounded  from above by $\delta(d).$
On the other hand, if $B$ is special and $\f O$ is an orbifold such  that 
$B:\f O\rightarrow \f O$ is a covering map between orbifolds, then $X:\f O\rightarrow \f O$ is also a covering map between orbifolds.
\et
\pr 
The first part of the theorem follows from Theorem \ref{main2} (see the introduction). So, assume that  
$B$ is special, and let $\f O=(\C\P^1\setminus\f E_B,\nu)$ be an orbifold such that $B:\f O\rightarrow \f O$ is a covering map.
We observe first that \be \l{epti} X^{-1}(\f E_B)=\f E_B.\ee
Indeed, since $B^{-1}(\f E_B)=\f E_B$, it follows from the commutativity of the diagram 
\be \l{ddii}
\begin{CD}
\C\P^1 @>B>> \C\P^1\\
@VV X V @VV X V\\ 
\C\P^1 @>B >> \C\P^1\ 
\end{CD}
\ee
that 
\be \l{epta2} B^{-1}(X^{-1}(\f E_B))=X^{-1}(\f E_B).\ee
Therefore, $X^{-1}(\f E_B)\subseteq \f E_B$. Since the set $X^{-1}(\f E_B)$ contains at least  $\vert \f E_B\vert$
points, this implies \eqref{epti}.

Set $R=\C\P^1\setminus \f E_B$ and $$F=B\circ X=X\circ B.$$
Since $X:R\rightarrow R$ and $F:R\rightarrow R$ are branched covering maps by \eqref{epti}, we can define
the orbifolds $X^*\f O$ and $F^*\f O.$  
By Theorem \ref{serrr}, we have:
$$F^*\f O=X^*(B^*\f O).$$   
However, since 
$B:\f O\rightarrow \f O$ is a covering map, the equality $B^*\f O=\f O$ holds, implying that $F^*\f O=X^*\f O.$
Thus, $F: X^*\f O\rightarrow \f O$ is a minimal holomorphic map between orbifolds. It follows now from $F=X\circ B$ by Corollary \ref{indu2} that  $B:X^*\f O\rightarrow X^*\f O$ also is  a 
minimal holomorphic map. In particular,  $\chi(X^*\f O)\geq 0$ by \eqref{iioopp}. 
On the other hand, since $X: X^*\f O\rightarrow \f O$ is a minimal holomorphic map, it follows from 
\eqref{iioopp} that  $\chi(X^*\f O)\leq 0.$  Therefore, \be \l{ik} \chi(X^*\f O)=0,\ee implying by    
Proposition \ref{p1} that $B:X^*\f O\rightarrow X^*\f O$ and $X: X^*\f O\rightarrow \f O$ are covering maps. 

If $\f E_B$ consists of two points, then  \eqref{ik} implies that the both orbifolds $\f O$ and $X^*\f O$ are non-ramified, so that 
\be \l{as} X^*\f O=\f O.\ee Thus, $X: \f O\rightarrow \f O$ is a covering map, as required. On the other hand, if  $\f E_B$ consists of one point, then 
without loss of generality we may assume that $\nu(-1)=\nu(1)=2$ and $B=\pm T_d,$ while \eqref{ik} implies that the orbifold 
$X^*\f O$ is defined by the equality $\hat \nu(a)=\hat \nu(b)=2$, where $a$ and $b$ are some points on $\C\P^1.$ It is not hard to see however that  $\pm T_d: X^*\f O\rightarrow X^*\f O$ cannot be a covering map for such $X^*\f O$  
unless $\{a,b\}=\{-1,1\}.$ Thus, in this case equality \eqref{as} is also satisfied. 

Finally, if $\f E_B=\emptyset,$ then $B$ is a Latt\'es map, and it is well-known that the orbifold $\f O$ such that $B:\f O\rightarrow \f O$ is a covering map is defined in a unique way by the dynamics of $B$ (see \cite{mil2}, or \cite{lattes}, Theorem 6.1). Thus, equality \eqref{as} still holds. \qed

\end{section}

\bibliographystyle{amsplain}

\end{document}